\documentclass[11pt]{amsart}
\usepackage{amssymb}

\usepackage{graphicx}
\usepackage{amscd}

\newtheorem{theorem}{Theorem}[section]
\theoremstyle{plain}
\newtheorem{lemma}[theorem]{Lemma}
\newtheorem{corollary}[theorem]{Corollary}
\newtheorem{proposition}[theorem]{Proposition}

\theoremstyle{definition}
\newtheorem{definition}[theorem]{Definition}
\newtheorem{example}[theorem]{Example}

\theoremstyle{remark}
\newtheorem{remark}[theorem]{Remark}
\numberwithin{equation}{section}

\input{tcilatex}

\begin{document}
\title[Nahm Algebras]{Nahm Algebras}
\author{Michael K. Kinyon}
\address{Department of Mathematics \& Computer Science\\
Indiana University South Bend\\
South Bend, IN 46634 USA}
\email{mkinyon@iusb.edu}
\urladdr{http://oit.iusb.edu/\symbol{126}mkinyon}
\author{Arthur A. Sagle}
\address{Department of Mathematics\\
University of Hawaii - Hilo\\
Hilo, HI 96720 USA}
\email{sagle@hawaii.edu}
\urladdr{}
\date{\today}
\subjclass{Primary: 17D, 17B; Secondary: 34A34}
\keywords{nonassociative algebras, Lie algebras, Nahm algebras.}
\thanks{The authors would like to thank N.C. Hopkins of Indiana State University for
some useful suggestions. One of us (AAS) would also like to thank the
University of Utah for its hospitality while this research was being carried
out.}

\begin{abstract}
Given a Lie algebra $\frak{g}$, the \emph{Nahm algebra} of $\frak{g}$ is the
vector space $\frak{g}\times \frak{g}\times \frak{g}$ with the natural
commutative, nonassociative algebra structure associated with the system of
ordinary differential equations (\ref{eqn1})-(\ref{eqn3}). Motivated by
applications to these equations, we herein initiate the study of Nahm
algebras.
\end{abstract}

\maketitle

\section{Introduction}

Let $\frak{g}$ be a real or complex Lie algebra. The \emph{Nahm equations}
for $\frak{g}$ are the following autonomous, first order differential
equations: 
\begin{eqnarray}
\dot{x} &=&[y,z]  \label{eqn1} \\
\dot{y} &=&[z,x]  \label{eqn2} \\
\dot{z} &=&[x,y]  \label{eqn3}
\end{eqnarray}
for $x,y,z\in \frak{g}$. This system of equations is of interest in
mathematical physics, especially in the case when $\frak{g}$ is a matrix Lie
algebra. This is because certain types of solutions of (\ref{eqn1})-(\ref
{eqn3}) are equivalent to monopole solutions of the self-dual $SU(2)$
Yang-Mills equations \cite{hitchin} \cite{nahm}. Much work has been done to
understand solutions of the Nahm equations in various physical and geometric
contexts; good places to start for those interested are the papers \cite
{biquard} \cite{hitchin} \cite{kovalev} \cite{nahm}, and the references
therein.

Let $Q:\frak{g}\times \frak{g}\times \frak{g}\rightarrow \frak{g}\times 
\frak{g}\times \frak{g}$ denote the mapping defined by the right hand side
of the Nahm equations (\ref{eqn1})-(\ref{eqn3}). This mapping is homogeneous
quadratic, i.e., $Q(\alpha X)=\alpha ^{2}Q(X)$ for all $\alpha \in \mathbb{K}
$, $X=(x_{1},x_{2},x_{3})\in \frak{g}\times \frak{g}\times \frak{g}$. In
1960, L. Markus \cite{markus} noted that to every quadratic differential
equation $\dot{X}=Q(X)$ occurring in a vector space $V$ over $\mathbb{K}$,
there is associated a natural algebra. This algebra is $\mathcal{A}=(V,\cdot
)$ where the operation $\cdot $ is the bilinearization of $Q$ defined by 
\begin{equation*}
X\cdot Y=\frac{1}{2}Q^{(2)}(0)(X,Y),
\end{equation*}
$X,Y\in V$, where $Q^{(2)}(0):V\times V\rightarrow V$ is the second
derivative of $Q$ at $0$. Clearly $\mathcal{A}$ is a commutative algebra,
but in general, it is nonassociative. We have $Q(X)=X\cdot X$, and if we
abbreviate $X^{2}:=X\cdot X$, then we may write the differential equation as 
$\dot{X}=X^{2}$. Thus we may view quadratic differential equations\ as
occurring in commutative, nonassociative algebras.

This algebraic perspective for quadratic differential equations is useful
because the structure of the underlying algebra can give information about
the trajectories (solution curves) of the differential equation. This is
analogous to the situation with constant coefficient linear differential
systems; such equations can be completely understood in terms of the theory
of a vector space acted on by a single linear transformation. It is
reasonable to expect that the theory of vector spaces with a bilinear
mapping, i.e., algebras, would play a role in understanding quadratic
differential equations. This line of investigation has been pursued by a
number of authors. For surveys, see \cite{ks1} or the monograph \cite{wbook}.

Applying these ideas to the Nahm equations (\ref{eqn1})-(\ref{eqn3}) leads
us to the following definition.

\begin{definition}
\label{defn-nahm-algebra}Let $(\frak{g},[\cdot ,\cdot ])$ be a real or
complex Lie algebra. The \emph{Nahm algebra} $(A(\frak{g}),\cdot )$
associated to $\frak{g}$ is the vector space $A(\frak{g})=\frak{g}\times 
\frak{g}\times \frak{g}$ with the multiplication 
\begin{equation}
X\cdot Y=\left( 
\begin{array}{c}
x_{1} \\ 
x_{2} \\ 
x_{3}
\end{array}
\right) \cdot \left( 
\begin{array}{c}
y_{1} \\ 
y_{2} \\ 
y_{3}
\end{array}
\right) \equiv \frac{1}{2}\left( 
\begin{array}{c}
\lbrack x_{2},y_{3}]+[y_{2},x_{3}] \\ 
\lbrack x_{3},y_{1}]+[y_{3},x_{1}] \\ 
\lbrack x_{1},y_{2}]+[y_{1},x_{2}]
\end{array}
\right)  \label{nahm-product}
\end{equation}
for $X=(x_{1},x_{2},x_{3})^{T},Y=(y_{1},y_{2},y_{3})^{T}\in \frak{g}\times 
\frak{g}\times \frak{g}$.
\end{definition}

In this paper, we initiate an investigation of Nahm algebras. As indicated,
our eventual goal is to understand the Nahm equations. Thus throughout the
paper, we will motivate the topics we discuss by referring to their
relevance for quadratic differential equations occurring in commutative
algebras. However, the paper itself is a purely algebraic study of Nahm
algebras. The implications of our results for the Nahm equations will appear
elsewhere.

Since our motivations lie in differential equations, all vector spaces and
algebras herein are over the field $\mathbb{K}$, where $\mathbb{K}=\mathbb{R}
$ or $\mathbb{K}=\mathbb{C}$. However, many of the results hold for
arbitrary fields of characteristic zero, and some hold in positive
characteristic. Definition \ref{defn-nahm-algebra} shows our notation
convention: lower case letters indicate elements of the Lie algebra $\frak{g}
$, and the corresponding upper case letters denote elements of the Nahm
algebra $A(\frak{g})$. We will frequently abbreviate the product in the Nahm
algebra by concatenation $XY=X\cdot Y$. We will also use the following
notation: for $i=1,2,3$, we define the projection $\pi _{i}:\frak{g}\times 
\frak{g}\times \frak{g}\rightarrow \frak{g}$ by $\pi
_{i}(x_{1},x_{2},x_{3})^{t}=x_{i}$.

If $\phi :\frak{g}\rightarrow \frak{h}$ is a homomorphism of Lie algebras,
then the mapping $A\left( \phi \right) :A\left( \frak{g}\right) \rightarrow
A\left( \frak{h}\right) $ defined by $A\left( \phi \right)
(x_{1},x_{2},x_{3})^{t}=(\phi \left( x_{1}\right) ,\phi \left( x_{2}\right)
,\phi \left( x_{3}\right) )^{t}$ is clearly a homomorphism of the associated
Nahm algebras. It follows that the assignment $\frak{g}\rightsquigarrow
A\left( \frak{g}\right) $ is a covariant functor from the category of Lie
algebras to the category of Nahm algebras.

One might guess that the Jacobi identity in the Lie algebra $\frak{g}$ would
lead to identities satisfied in the Nahm algebra $A(\frak{g})$.
Interestingly enough, this does not seem to be the case. For instance, the
Nahm product (\ref{nahm-product}) is not, in general,\ fourth
power-associative, and thus Nahm algebras are not a subclass of some
well-studied variety of commutative, power-associative algebras, such as
Jordan algebras \cite{schafer}.

If $T:\mathbb{K}^{3}\rightarrow \mathbb{K}^{3}$ is a linear transformation,
then $T$ acts on $A(\frak{g})=\frak{g}\otimes \mathbb{K}^{3}$ in the obvious
way: 
\begin{equation*}
T\left( 
\begin{array}{c}
x_{1} \\ 
x_{2} \\ 
x_{3}
\end{array}
\right) =x_{1}\otimes T\mathbf{e}_{1}+x_{2}\otimes T\mathbf{e}%
_{2}+x_{3}\otimes T\mathbf{e}_{3}
\end{equation*}
where $\mathbf{e}_{i}$ is the $i$th standard basis vector of $\mathbb{K}^{3}$%
. More specifically, if $T$ is given by a $3\times 3$ matrix $T=[T_{ij}]$
relative to the standard basis, then the action agrees with that obtained by
formally multiplying the column vector $(x_{1},x_{2},x_{3})^{t}$ ($x_{i}\in 
\frak{g}$, $i=1,2,3$) on the left by the matrix $T$.

For any linear transformation $L\in gl(\frak{g})$, we will denote the
naturally induced transformation in $gl(A(\frak{g}))$ by $\mathrm{diag}(L)$;
thus 
\begin{equation*}
\mathrm{diag}(L)\left( 
\begin{array}{c}
x_{1} \\ 
x_{2} \\ 
x_{3}
\end{array}
\right) =\left( 
\begin{array}{c}
Lx_{1} \\ 
Lx_{2} \\ 
Lx_{3}
\end{array}
\right)
\end{equation*}
for $X=(x_{1},x_{2},x_{3})^{T}\in A(\frak{g})$.

For $X=(x_{1},x_{2},x_{3})^{T}\in A(\frak{g})$, the \emph{left
multiplication operator }$L(X)\in gl(A(\frak{g}))$ is defined by $L(X)Y=XY$
for all $Y\in A(\frak{g})$. Using (\ref{nahm-product}), we see that $L(X)$
has a block matrix representation given by 
\begin{equation}
L(X)=\frac{1}{2}\left[ 
\begin{array}{ccc}
0 & -\mathrm{ad}\ x_{3} & \mathrm{ad}\ x_{2} \\ 
\mathrm{ad}\ x_{3} & 0 & -\mathrm{ad}\ x_{1} \\ 
-\mathrm{ad}\ x_{2} & \mathrm{ad}\ x_{1} & 0
\end{array}
\right]  \label{L(X)}
\end{equation}
where for $x\in \frak{g}$, the adjoint representation is given by $(\mathrm{%
ad}\ x)y=[x,y]$.

This suggests the following definition. Let $\rho :\frak{g}\rightarrow gl(V)$
be a representation of $\frak{g}$ on a vector space $V$. For $X\in A(\frak{g}%
)$, we define an operator $L_{\rho }(X)\in gl(V\times V\times V)$ as
follows: 
\begin{equation}
L_{\rho }(X)=\frac{1}{2}\left[ 
\begin{array}{ccc}
0 & -\rho (x_{3}) & \rho (x_{2}) \\ 
\rho \left( x_{3}\right) & 0 & -\rho \left( x_{1}\right) \\ 
-\rho \left( x_{2}\right) & \rho \left( x_{1}\right) & 0
\end{array}
\right] \text{.}  \label{L-rho}
\end{equation}
Thus $L_{\mathrm{ad}}(X)=L(X)$. We will use the operators $L_{\rho }(X)$, $%
X\in A(\frak{g})$, in our discussion of invariant bilinear forms.

We conclude this introduction with an outline of the sequel. In \S 2 we
discuss the relationship between subalgebras and ideals of the Lie algebra $%
\frak{g}$ and subalgebras and ideals of the Nahm algebra $A(\frak{g})$. In
\S 3 we discuss how $\mathbb{Z}_{2}$-gradings of $\frak{g}$ naturally induce 
$\mathbb{Z}_{2}$-gradings of $A(\frak{g})$. We also show that every Nahm
algebra $A(\frak{g})$ has a natural $\mathbb{Z}_{2}$-grading where the even
subalgebra is a copy (as a vector space) of $\frak{g}$ itself. In \S 4 we
discuss nilpotents of index 2 and idempotents in $A(\frak{g})$. Roughly
speaking, nilpotents in $A(\frak{g})$ are built from abelian subalgebras of $%
\frak{g}$, and idempotents in $A(\frak{g})$ are built from subalgebras of $%
\frak{g}$ which are isomorphic to $so(3,\mathbb{K})$. In \S 5, we discuss
simplicity and prove that the Lie algebra $\frak{g}$ is simple if and only
if the Nahm algebra $A(\frak{g})$ is simple. We also give an example to show
that simple Nahm algebras can have simple subalgebras which are not
themselves Nahm algebras of a Lie algebra. We then turn to semisimplicity
and prove that the Lie algebra $\frak{g}$ is semisimple if and only if the
Nahm algebra $A(\frak{g})$ is semisimple. In \S 6 we show that the radical
of a Nahm algebra is the Nahm algebra of the radical of the Lie algebra. It
follows from this that every Nahm algebra has a Levi-Malcev decomposition.
In \S 7 we consider invariant bilinear forms for Nahm algebras. Any
representation of $\frak{g}$ naturally induces an invariant trace form on $A(%
\frak{g})$. The form so induced by the adjoint representation, which we call
the standard form, measures the semisimplicity of $A(\frak{g})$ in exact
analogy with the role of the Killing form on $\frak{g}$ itself: $A(\frak{g})$
is semisimple if and only if its standard form is nondegenerate. In \S 8 we
consider derivations of $A(\frak{g})$. We show that the derivation algebra
of any Nahm algebra has two natural subalgebras: one is a copy of \textrm{ad}%
$(\frak{g})$, and the other is a copy of $so(3,\mathbb{K})$. For $A(\frak{g}%
) $ simple, we prove that the derivation algebra is exactly the direct sum
of these two subalgebras. We first show this for $\mathbb{K}=\mathbb{C}$ and
then note that the result follows in the real case by the invariance of
dimension of derivation algebras. Along the way, we also prove a version of
Schur's lemma for complex, simple Nahm algebras. Finally, in \S 9 we discuss
automorphisms of Nahm algebras, and we characterize the automorphism group
of a simple Nahm algebra: it is a direct product of the automorphism group
of the Lie algebra and $SO(3,\mathbb{K})$.

We should mention that for any anticommutative algebra $\frak{g}$, one can
certainly define its ``Nahm algebra'' $A(\frak{g})$ as in Definition \ref
{defn-nahm-algebra}. Indeed, some of what follows is valid in the case
where, for example, $\frak{g}$ is a Malcev algebra. However, in this paper $%
\frak{g}$ will be a Lie algebra, and the structure of the Nahm algebra will
turn out to be closely related to the structure of $\frak{g}$ itself.

\section{Subalgebras and Ideals}

A \emph{subalgebra} of an algebra $\mathcal{A}$ is a subspace $\mathcal{B}$
such that $\mathcal{B}^{2}\subseteq \mathcal{B}$, that is, $XY\in \mathcal{B}
$ for all $X,Y\in \mathcal{B}$. Given a fixed $P\in \mathcal{A}$, the \emph{%
subalgebra generated by }$P$ is defined by 
\begin{equation}
\langle P\rangle =\mathbb{R}\left[ P\right] ,  \label{subalg(P)}
\end{equation}
which is the set of all polynomials in $P$. For a commutative algebra $%
\mathcal{A}$ with its associated quadratic differential equation $\dot{X}%
=X^{2}$, the unique solution $X\left( t;P\right) $ to the equation
satisfying the initial condition $X\left( 0\right) =P$ remains in the
subalgebra $\langle P\rangle $, i.e., $X\left( t;P\right) \in \langle
P\rangle $ for all $t$ \cite{markus} \cite{ks1} \cite{wbook}.

For a Nahm algebra $A=A(\frak{g})$ of a Lie algebra $\frak{g}$, it is
reasonable to expect that subalgebras in $A$ are related to subalgebras in $%
\frak{g}$.

\begin{theorem}
\label{thm-subalg}Let $\frak{m}_{i}\subseteq \frak{g}$, $i=1,2,3$, be
subspaces, and let $M=\frak{m}_{1}\times \frak{m}_{2}\times \frak{m}_{3}$.
Then $M$ is a subalgebra $A\left( \frak{g}\right) $ if and only if $[\frak{m}%
_{i},\frak{m}_{i+1}]\subseteq \frak{m}_{i+2}$ for $i=1,2,3$ (where index
addition is modulo $3$).
\end{theorem}

\begin{proof}
This follows immediately from considering the components of the product $XY$
for $X,Y\in B$.
\end{proof}

\begin{remark}
\label{rem-subalg}It is important to note that Theorem \ref{thm-subalg} does
not characterize arbitrary subalgebras of $A\left( \frak{g}\right) $; it
only characterizes those which are direct products of subspaces in $\frak{g}$%
.
\end{remark}

\begin{corollary}
\label{coro-subalg}Let $\frak{m}\subseteq \frak{g}$ be a subspace. Then $M=%
\frak{m}\times \frak{m}\times \frak{m}$ is a subalgebra of $A$ if and only
if $\frak{m}$ is a subalgebra of $\frak{g}$. In this case, $M=A\left( \frak{m%
}\right) $ is the Nahm algebra of $\frak{m}$.
\end{corollary}

An \emph{ideal }of an algebra $\mathcal{A}$ is a subspace $\mathcal{J}$ such
that $\mathcal{J}\mathcal{A}=\mathcal{AJ}\subseteq \mathcal{J}$, that is, $%
XY\in \mathcal{J}$ and $YX\in \mathcal{J}$ for all $X\in \mathcal{A},Y\in 
\mathcal{J}$. (Since we are dealing only with commutative and Lie algebras
here, ``ideal'' for us means ``two-sided ideal''.) In case $\mathcal{A}$ is
a commutative algebra, the presence of an ideal $\mathcal{J}$ in $\mathcal{A}
$ implies that the associated quadratic differential equation can be
decoupled into a quadratic equation in the quotient space $\mathcal{A}/%
\mathcal{J}$ and a (nonhomogeneous) quadratic differential equation in $%
\mathcal{J}$; see \cite{wbook}, p.23.

Let $A(\frak{g})$ be the Nahm algebra of a Lie algebra $\frak{g}$. As one
might expect, ideals in $A(\frak{g})$ are closely related to ideals in $%
\frak{g}$.

\begin{theorem}
\label{thm-ideals}Let $J$ be an ideal of $A(\frak{g})$, and let $\frak{h}%
_{i}=\pi _{i}(J)\subseteq \frak{g}$. Then $[\frak{g},\frak{h}_{i}]\subseteq 
\frak{h}_{i+1}\cap \frak{h}_{i+2}$ for $i=1,2,3$ (where index addition is
modulo $3$).
\end{theorem}

\begin{proof}
Fix $y_{1}\in \frak{h}_{1}$. There exists $y_{j}\in \frak{h}_{j}$, $j=2,3$,
such that $(y_{1},y_{2},y_{3})^{t}\in J$. For all $x\in $ $\frak{g}$, 
\begin{equation*}
\left( 
\begin{array}{c}
0 \\ 
x \\ 
x
\end{array}
\right) \cdot \left( 
\begin{array}{c}
y_{1} \\ 
y_{2} \\ 
y_{3}
\end{array}
\right) =\frac{1}{2}\left( 
\begin{array}{c}
\lbrack x,y_{3}]+[y_{2},x] \\ 
\lbrack x,y_{1}] \\ 
\lbrack y_{1},x]
\end{array}
\right) \in J\text{.}
\end{equation*}
Thus $[x,y_{1}]\in \frak{h}_{2}\cap \frak{h}_{3}$. Similar calculations show
the other inclusions.
\end{proof}

\begin{corollary}
\label{coro-intersect}Let $J$ be an ideal of $A(\frak{g})$ and let $\frak{h}%
_{i}=\pi _{i}(J)$, $i=1,2,3$. then $\frak{h}_{1}\cap \frak{h}_{2}\cap \frak{h%
}_{3}$ is an ideal of $\frak{g}$.
\end{corollary}

\begin{theorem}
\label{thm-ideals2}Let $\frak{h}_{i}\subseteq \frak{g}$, $i=1,2,3$, be
subspaces, and let $J=\frak{h}_{1}\times \frak{h}_{2}\times \frak{h}_{3}$.
Then $J$ is an ideal of $A(\frak{g})$ if and only if $[\frak{g},\frak{h}%
_{i}]\subseteq \frak{h}_{i+1}\cap \frak{h}_{i+2}$ for $i=1,2,3$ (where index
addition is modulo $3$).
\end{theorem}

\begin{proof}
The necessity of the stated condition is Theorem \ref{thm-ideals}, while the
sufficiency is clear from (\ref{nahm-product}).
\end{proof}

\begin{corollary}
\label{coro-ideals}Let $\frak{h}\subseteq \frak{g}$ be a subspace. Then $%
\frak{h}\times \frak{h}\times \frak{h}$ is an ideal of $A(\frak{g})$ if and
only if $\frak{h}$ is an ideal of $\frak{g}$.
\end{corollary}

\section{$\mathbb{Z}_{2}$-Gradings}

When a Lie algebra $\frak{g}$ has a $\mathbb{Z}_{2}$-grading, it induces an
interesting class of subalgebras of it associated Nahm algebra $A(\frak{g})$%
. Thus assume $\frak{g}=\frak{g}_{0}\oplus \frak{g}_{1}$ where $\frak{g}_{0}$
is a subalgebra, $\frak{g}_{1}$ is a subspace, $[\frak{g}_{0},\frak{g}%
_{1}]\subseteq \frak{g}_{1}$ and $[\frak{g}_{1},\frak{g}_{1}]\subseteq \frak{%
g}_{0}$. (For example, $\frak{g}$ could be a semisimple Lie algebra with
Cartan decomposition $\frak{g}=\frak{g}_{0}\oplus \frak{g}_{1}$.) Let 
\begin{equation}
A_{011}=\frak{g}_{0}\times \frak{g}_{1}\times \frak{g}_{1}  \label{A011}
\end{equation}
and similarly define $A_{101}$ and $A_{110}$. By Theorem \ref{thm-subalg}, $%
A_{011}$, $A_{101}$ and $A_{110}$ are subalgebras of $A(\frak{g})$. Let 
\begin{equation}
A_{100}=\frak{g}_{1}\times \frak{g}_{0}\times \frak{g}_{0}  \label{A100}
\end{equation}
and similarly define $A_{010}$ and $A_{001}$. Then the following properties
are easily seen to hold: 
\begin{eqnarray}
A(\frak{g}) &=&A_{011}\oplus A_{100}  \label{grade1} \\
A_{011}\cdot A_{100} &\subseteq &A_{100}  \label{grade2} \\
A_{100}\cdot A_{100} &\subseteq &A_{011}.  \label{grade3}
\end{eqnarray}
Similar results hold for the other subalgebras and subspaces. Therefore the $%
\mathbb{Z}_{2}$-grading $\frak{g}=\frak{g}_{0}\oplus \frak{g}_{1}$ naturally
induces three $\mathbb{Z}_{2}$-gradings of the Nahm algebra $A(\frak{g})$.
If a commutative algebra $\mathcal{A}$ has a $\mathbb{Z}_{2}$-grading $%
\mathcal{A}=\mathcal{A}_{0}\oplus \mathcal{A}_{1}$, then the differential
equation $\dot{X}=X^{2}$ in $\mathcal{A}$ can be decomposed relative to the
grading, and this can give information about the trajectories \cite{hk2}.

Every Nahm algebra carries a natural $\mathbb{Z}_{2}$-grading whether the
underlying Lie algebra is $\mathbb{Z}_{2}$-graded or not. For each $x\in 
\frak{g}$, let 
\begin{equation}
\Delta (x)=\left( 
\begin{array}{c}
x \\ 
x \\ 
x
\end{array}
\right)  \label{D(x)}
\end{equation}
and let 
\begin{equation}
\Delta \left( \frak{g}\right) =\left\{ \Delta (x):x\in \frak{g}\right\} .
\label{D(g)}
\end{equation}
As a subspace, $\Delta \left( \frak{g}\right) $\ is just a copy of $\frak{g}$
itself. Let 
\begin{equation*}
W\left( \frak{g}\right) =\left\{ X\in A\left( \frak{g}\right)
:x_{1}+x_{2}+x_{3}=0\right\} .
\end{equation*}
Then the following properties hold.

\begin{theorem}
\label{thm-A(g)=D(g)(+)W(g)}

\begin{enumerate}
\item  $\Delta \left( \frak{g}\right) $ is an abelian subalgebra of $A\left( 
\frak{g}\right) $, i.e., 
\begin{equation*}
\Delta (x)\Delta (y)=0
\end{equation*}
for all $x,y\in \frak{g}$.

\item  $A\left( \frak{g}\right) =\Delta \left( \frak{g}\right) \oplus
W\left( \frak{g}\right) .$

\item  $\Delta \left( \frak{g}\right) \cdot W\left( \frak{g}\right)
\subseteq W\left( \frak{g}\right) .$

\item  $W\left( \frak{g}\right) \cdot W\left( \frak{g}\right) \subseteq
\Delta \left( \frak{g}\right) .$
\end{enumerate}
\end{theorem}

\begin{proof}
1. This is clear from (\ref{nahm-product}).

2. If $\Delta (x)\in \Delta \left( \frak{g}\right) \cap W\left( \frak{g}%
\right) $, then $3x=0$ and thus $\Delta (x)=0$. For $X\in A\left( \frak{g}%
\right) $, define mappings $P_{\Delta }:A\left( \frak{g}\right) \rightarrow
\Delta \left( \frak{g}\right) $ and $P_{W}:A\left( \frak{g}\right)
\rightarrow W\left( \frak{g}\right) $ by 
\begin{equation}
P_{\Delta }(X)=\Delta \left( \frac{1}{3}(x_{1}+x_{2}+x_{3})\right)
\label{P_D}
\end{equation}
and 
\begin{equation}
P_{W}(X)=\frac{1}{3}\left( 
\begin{array}{c}
2x_{1}-x_{2}-x_{3} \\ 
-x_{1}+2x_{2}-x_{3} \\ 
-x_{1}-x_{2}+2x_{3}
\end{array}
\right) \text{.}  \label{P_W}
\end{equation}
Then $P_{\Delta }^{2}=P_{\Delta }$, $P_{W}^{2}=P_{W}$, and $X=P_{\Delta
}(X)+P_{W}(X)$ for each $X\in A\left( \frak{g}\right) $, which proves the
desired result and also shows that $P_{\Delta }$ and $P_{W}$ are the
projectors onto $\Delta \left( \frak{g}\right) $ and $W\left( \frak{g}%
\right) $, respectively.

3. For $x\in \frak{g}$, $Y\in W\left( \frak{g}\right) $, we have 
\begin{equation*}
\Delta (x)Y=\frac{1}{2}\left( 
\begin{array}{c}
\lbrack x,y_{3}-y_{2}] \\ 
\lbrack x,y_{1}-y_{3}] \\ 
\lbrack x,y_{2}-y_{1}]
\end{array}
\right) \text{.}
\end{equation*}
Summing up entries shows $\Delta (x)Y\in W\left( \frak{g}\right) $.

4. Fix $X,Y\in W\left( \frak{g}\right) $, and let $s_{12}$ denote the
difference between the first and second entries of the product $2XY$. Then 
\begin{eqnarray*}
s_{12} &=&\left( [x_{2},y_{3}]+[y_{2},x_{3}]\right) -\left(
[x_{3},y_{1}]+[y_{3},x_{1}]\right) \\
&=&[x_{1}+x_{2},y_{3}]+[y_{1}+y_{2},x_{3}] \\
&=&[x_{1}+x_{2},y_{3}]+[y_{1}+y_{2},x_{3}],
\end{eqnarray*}
using anticommutativity of the Lie bracket. Adding the quantity $%
[x_{3},y_{3}]+[y_{3},x_{3}]=0$, we obtain 
\begin{eqnarray*}
s_{12} &=&[x_{1}+x_{2}+x_{3},y_{3}]+[y_{1}+y_{2}+y_{3},x_{3}] \\
&=&0,
\end{eqnarray*}
since $X,Y\in W\left( \frak{g}\right) $. Similarly, we can show that the
difference between the second and third entries of $2XY$ is $0$. Thus all
three entries are equal, which proves that $XY\in \Delta \left( \frak{g}%
\right) $.
\end{proof}

\begin{corollary}
\label{coro-grading}The decomposition $A\left( \frak{g}\right) =\Delta
\left( \frak{g}\right) \oplus W\left( \frak{g}\right) $ is a $\mathbb{Z}_{2}$%
-grading for $A\left( \frak{g}\right) $ with even subalgebra $\Delta \left( 
\frak{g}\right) $ and odd subspace $W\left( \frak{g}\right) $.
\end{corollary}

\begin{definition}
\label{defn-diagonal}$\Delta \left( \frak{g}\right) \subset A\left( \frak{g}%
\right) $ is called the \emph{diagonal subalgebra} of $A\left( \frak{g}%
\right) $.
\end{definition}

The decomposition $A\left( \frak{g}\right) =\Delta \left( \frak{g}\right)
\oplus W\left( \frak{g}\right) $ of a Nahm algebra has implications for the
Nahm equations as we will discuss elsewhere.

\section{Idempotents and Nilpotents}

A \emph{nilpotent} element $N$ (of index $2$) of an algebra $\mathcal{A}$ is
one which satisfies $N^{2}=0$. If $\mathcal{A}$ is the commutative algebra
associated to a quadratic differential equation $\dot{X}=X^{2}$, then
nilpotents correspond to nonzero equilibria, i.e., stationary points. In
particular, if $N$ is a nilpotent, then the constant function $X(t)=N$ is a
solution: 
\begin{equation*}
N^{2}=X^{2}=\frac{dX}{dt}=\frac{dN}{dt}=0.
\end{equation*}
Conversely, the same calculation shows that a constant function $X\left(
t\right) =N$ is a solution only if $N^{2}=0$.

Now let $\frak{g}$ be a Lie algebra with Nahm algebra $A=A\left( \frak{g}%
\right) $. Assume that $N=(n_{1},n_{2},n_{3})^{t}\in A$ is a nilpotent.
Since $N^{2}=0$, we have $[n_{i},n_{j}]=0$ for all $i,j$. Thus the Lie
subalgebra $\langle n_{1},n_{2},n_{3}\rangle $ of $\frak{g}$ generated by
the set $\left\{ n_{1},n_{2},n_{3}\right\} $ is abelian, and the subspace $%
\mathbb{K}\cdot n_{1}\times \mathbb{K}\cdot n_{2}\times \mathbb{K}\cdot
n_{3}\subseteq A$ is an abelian subalgebra. Conversely, given any abelian
subalgebra $\frak{h}$ of $\frak{g}$, any three-element set $\left\{
n_{1},n_{2},n_{3}\right\} \subset \frak{h}$ gives an abelian subalgebra $%
\mathbb{K}\cdot n_{1}\times \mathbb{K}\cdot n_{2}\times \mathbb{K}\cdot
n_{3} $ of $A$, every element of which is a nilpotent. Indeed, if $B$ is an
abelian subalgebra of $A$, then clearly every element of $B$ is a nilpotent.
This applies, for instance, to the diagonal subalgebra $\Delta \left( \frak{g%
}\right) $.

An \emph{idempotent} of an algebra $\mathcal{A}$ is a nonzero element $E\in 
\mathcal{A}$ satisfying $E^{2}=E$. If $\mathcal{A}$ is the commutative
algebra associated to a quadratic differential equation $\dot{X}=X^{2}$,
then idempotents gives solutions which blow up in finite time along rays.
For arbitrary $E\in \mathcal{A}$ and for $a\neq 0$, the function $%
X(t)=aE/(1-at)$ blows up in finite time at $t=1/a$. If $E$ is an idempotent,
then $X\left( t\right) $ is the solution to the differential equation with
initial value $aE$: 
\begin{equation*}
\frac{a^{2}E^{2}}{\left( 1-at\right) ^{2}}=X^{2}=\frac{dX}{dt}=\frac{a^{2}E}{%
\left( 1-at\right) ^{2}}.
\end{equation*}
Conversely, the same calculation shows that $X\left( t\right) $ is a
solution only if $E=E^{2}$ is an idempotent.

Let $\frak{g}$ be a Lie algebra with Nahm algebra $A=A\left( \frak{g}\right) 
$. Assume that $E=(e_{1},e_{2},e_{3})^{t}\in A\left( \frak{g}\right) $ is an
idempotent. Since $E=E^{2}$, we have 
\begin{equation}
\left( 
\begin{array}{c}
e_{1} \\ 
e_{2} \\ 
e_{3}
\end{array}
\right) =\left( 
\begin{array}{c}
\lbrack e_{2},e_{3}] \\ 
\lbrack e_{3},e_{1}] \\ 
\lbrack e_{1},e_{2}]
\end{array}
\right) \text{.}  \label{su(2)}
\end{equation}
If any $e_{i}=0$, then (\ref{su(2)}) shows that all $e_{i}=0$, which
contradicts $E$ being an idempotent. Suppose $\left\{
e_{1},e_{2},e_{3}\right\} $ is linearly dependent, say, $e_{3}=ae_{1}+be_{2}$
for some $a,b\in \mathbb{K}$. Then 
\begin{equation*}
e_{1}=[e_{2},ae_{1}+be_{2}]=-a[e_{1},e_{2}]=-ae_{3}=-a^{2}e_{1}-abe_{2},
\end{equation*}
and 
\begin{equation*}
e_{2}=[ae_{1}+be_{2},e_{1}]=-b[e_{1},e_{2}]=-be_{3}=-abe_{1}-b^{2}e_{2}.
\end{equation*}
Thus $a^{2}=b^{2}=-1$ and $ab=0$, which is a contradiction. Therefore $%
\left\{ e_{1},e_{2},e_{3}\right\} $ is linearly independent, and (\ref{su(2)}%
) shows that it satisfies the defining relations of the Lie algebra $so(3,%
\mathbb{K})$ (or isomorphically, the Lie algebra $\mathbb{K}^{3}$ with the
cross-product as Lie bracket).

Conversely, assume that $\frak{g}$ contains a subalgebra $B$ isomorphic to $%
so(3,\mathbb{K})$. Let $\left\{ e_{1},e_{2},e_{3}\right\} $ be an ordered
basis for $B$ satisfying the relations $[e_{i},e_{i+1}]=e_{i+2}$ (where
index addition is modulo $3$). Set $E=(e_{1},e_{2},e_{3})^{t}$. Then $E$ is
an idempotent.

\section{Simple and Semisimple Algebras}

An algebra $\mathcal{A}$ is called \emph{simple} if $\mathcal{A}^{2}\neq
\left\{ 0\right\} $ and $\mathcal{A}$ contains no nontrivial ideals, that
is, the only nonzero ideal of $\mathcal{A}$ is $\mathcal{A}$ itself. We now
consider the relationship between simplicity of a Lie algebra $\frak{g}$ and
the simplicity of its Nahm algebra $A(\frak{g})$.

\begin{theorem}
\label{thm-simple}$A(\frak{g})$ is simple if and only if $\frak{g}$ is
simple.
\end{theorem}

\begin{proof}
Assume $A(\frak{g})$ is simple and let $\frak{h}$ be a nonzero ideal of $%
\frak{g}$. By Corollary \ref{coro-ideals}, $\frak{h}\times \frak{h}\times 
\frak{h}$ is an ideal of $A(\frak{g})$. Thus $A(\frak{g})=\frak{h}\times 
\frak{h}\times \frak{h}$, which implies $\frak{h}=\frak{g}$. Conversely,
assume $\frak{g}$ is simple, let $J$ be an ideal of $A(\frak{g})$, and let $%
\frak{h}_{i}=\pi _{i}(J)$. By Theorem \ref{thm-ideals}, $[\frak{g},[\frak{g},%
\frak{h}_{i}]]\subseteq \frak{h}_{1}\cap \frak{h}_{2}\cap \frak{h}_{3}$ for $%
i=1,2,3$. By Corollary \ref{coro-intersect}, $\frak{h}_{1}\cap \frak{h}%
_{2}\cap \frak{h}_{3}$ is an ideal. If $\frak{h}_{1}\cap \frak{h}_{2}\cap 
\frak{h}_{3}=\frak{g}$, then $J=\frak{g}\times \frak{g}\times \frak{g}=A(%
\frak{g})$, and thus we may assume $\frak{h}_{1}\cap \frak{h}_{2}\cap \frak{h%
}_{3}=\left\{ 0\right\} $ since $\frak{g}$ is simple. Now for $x,y\in \frak{g%
}$, $z\in \frak{h}_{i}$, the Jacobi identity gives $%
[[x,y],z]=[[x,z],y]+[y,[x,z]]\in $ $[\frak{g},[\frak{g},\frak{h}%
_{i}]]=\left\{ 0\right\} $. Therefore $[[\frak{g},\frak{g}],\frak{h}%
_{i}]]=\left\{ 0\right\} $. But $\frak{g=}[\frak{g},\frak{g}]$ since $\frak{g%
}$ is simple, and thus $[\frak{g},\frak{h}_{i}]=\left\{ 0\right\} $. Hence
each $\frak{h}_{i}$ is contained in the center of $\frak{g}$, and thus each $%
\frak{h}_{i}=\left\{ 0\right\} $. Therefore $J=\left\{ 0\right\} $, which
shows that $A(\frak{g})$ is simple.
\end{proof}

The next example shows that simple Nahm algebras can have simple subalgebras
which are not themselves Nahm algebras of a Lie algebra.

\begin{example}
Let $\frak{g}=\mathbb{R}^{3}$ with the bracket being the cross-product.
(Thus $\frak{g}$ is isomorphic to $so(3)$.) Then $\frak{g}$ is simple, and
thus the Nahm algebra $A\left( \frak{g}\right) =\mathbb{R}^{3}\times \mathbb{%
R}^{3}\times \mathbb{R}^{3}$ is simple by Theorem \ref{thm-simple}. Now let $%
\mathbf{e}_{1},\mathbf{e}_{2},\mathbf{e}_{3}$ denote the standard basis of $%
\mathbb{R}^{3}$, and consider the subspace $B=\mathbb{R}\cdot \mathbf{e}%
_{1}\times \mathbb{R}\cdot \mathbf{e}_{2}\times \mathbb{R}\cdot \mathbf{e}%
_{3}\subseteq A\left( \frak{g}\right) $. Then $B$ is clearly a subalgebra,
but it is not the Nahm algebra of a subalgebra of $\frak{g}$ (or of any Lie
algebra, for that matter). In addition, $B$ is simple as the following
argument shows. Let $J\subseteq B$ be a nonzero ideal, and assume $\mathbf{0}%
\neq (a\mathbf{e}_{1},b\mathbf{e}_{2},c\mathbf{e}_{3})^{t}\in J$. Assume
first that $a\neq 0$, $b=c=0$. Then 
\begin{equation*}
(\mathbf{0},\mathbf{e}_{2},\mathbf{0})^{t}\cdot (a\mathbf{e}_{1},\mathbf{0},%
\mathbf{0})^{t}=\frac{1}{2}(\mathbf{0},\mathbf{0},a\mathbf{e}_{3})^{t}\in J
\end{equation*}
and 
\begin{equation*}
(\mathbf{0},\mathbf{0},\mathbf{e}_{3})^{t}\cdot (a\mathbf{e}_{1},\mathbf{0},%
\mathbf{0})^{t}=\frac{1}{2}(\mathbf{0},a\mathbf{e}_{2},\mathbf{0})^{t}\in J.
\end{equation*}
This shows $J=B$. Applying similar arguments shows that we may assume that
at least two of $a,b,c$ are nonzero. Thus assume $a=0,b,c\neq 0$. Then 
\begin{equation*}
(\mathbf{0},\mathbf{e}_{2},\mathbf{0})^{t}\cdot (\mathbf{0},b\mathbf{e}_{2},c%
\mathbf{e}_{3})^{t}=\frac{1}{2}(c\mathbf{e}_{1},\mathbf{0},\mathbf{0}%
)^{t}\in J,
\end{equation*}
and repeating the argument above gives $J=B$. Applying similar arguments
shows that we may assume that each of $a,b,c$ are nonzero. Then 
\begin{equation*}
(\mathbf{0},b\mathbf{e}_{2},-c\mathbf{e}_{3})^{t}\cdot (a\mathbf{e}_{1},b%
\mathbf{e}_{2},c\mathbf{e}_{3})^{t}=\frac{1}{2}(\mathbf{0},-ca\mathbf{e}%
_{2},ab\mathbf{e}_{3})^{t}\in J,
\end{equation*}
and repeating the preceding arguments gives $J=B$. Thus $B$ is simple, as
claimed.
\end{example}

An algebra $\mathcal{A}$ is called \emph{semisimple} if there exist ideals $%
\mathcal{A}_{1},\ldots ,\mathcal{A}_{n}$, each of which is a simple algebra,
such that 
\begin{equation}
\mathcal{A}=\mathcal{A}_{1}\oplus \cdots \oplus \mathcal{A}_{n},
\label{semisimple}
\end{equation}
a direct sum of subalgebras. Any ideal $\mathcal{J}$ of $\mathcal{A}$ is
given by a direct sum $\mathcal{J}=\mathcal{A}_{i_{1}}\oplus \cdots \oplus 
\mathcal{A}_{i_{k}}$ for suitable $\mathcal{A}_{i_{j}}$. In particular, 
\begin{equation}
\mathcal{J}^{2}=\mathcal{A}_{i_{1}}^{2}\oplus \cdots \oplus \mathcal{A}%
_{i_{k}}^{2}=\mathcal{A}_{i_{1}}\oplus \cdots \oplus \mathcal{A}_{i_{k}}=%
\mathcal{J}  \label{Jsqrd=J}
\end{equation}
using $\mathcal{A}_{p}\mathcal{A}_{q}=0$ if $p\neq q$ and $\mathcal{A}%
_{p}^{2}=\mathcal{A}_{p}$.

If $\mathcal{A}$ is a semisimple commutative algebra, then the associated
quadratic differential equation $\dot{X}=X^{2}$ in $\mathcal{A}$ can be
decoupled into differential equations occurring in the simple subalgebras $%
\mathcal{A}_{j}$. This follows from the remarks above about ideals (see \cite
{wbook}, p.23), but can be just as easily seen directly. Thus if $X(t)$ is a
solution and $X=X_{1}+\cdots +X_{n}$ is the decomposition of $X$ given by (%
\ref{semisimple}), then 
\begin{equation*}
\dot{X}_{1}+\cdots +\dot{X}_{n}=\dot{X}=X^{2}=(X_{1}+\cdots
+X_{n})^{2}=X_{1}^{2}+\cdots +X_{n}^{2}\text{,}
\end{equation*}
since $X_{i}X_{j}=0$ for $i\neq j$. Thus $\dot{X}_{j}=X_{j}^{2}$ for $%
j=1,\ldots ,n$. See \cite{ks1} \cite{wbook}.

In analogy with our discussion of simplicity, we now show the relationship
between the semisimplicity of a Lie algebra $\frak{g}$ and the
semisimplicity of its associated Nahm algebra $A\left( \frak{g}\right) $.

\begin{theorem}
\label{thm-semisimple}$A(\frak{g})$ is semisimple if and only if $\frak{g}$
is semisimple.
\end{theorem}

\begin{proof}
Let $\frak{g}=\frak{g}_{1}\oplus \cdots \oplus \frak{g}_{n}$ be semisimple
with each $\frak{g}_{i}$ simple. By Corollary \ref{coro-ideals} and Theorem 
\ref{thm-simple}, each $\frak{g}_{i}\times \frak{g}_{i}\times \frak{g}_{i}$
is an ideal of $A(\frak{g})$ and a simple algebra. Since $A(\frak{g}%
)=\bigoplus_{i=1}^{n}\left( \frak{g}_{i}\times \frak{g}_{i}\times \frak{g}%
_{i}\right) $, a direct sum, $A(\frak{g})$ is semisimple.

Conversely, assume $A(\frak{g})=A_{1}\oplus \cdots \oplus A_{n}$ is
semisimple with each $A_{i}$ simple. Let $\frak{h}\neq \left\{ 0\right\} $
be a solvable ideal of $\frak{g}$. Then 
\begin{equation*}
\frak{h}\supset \frak{h}^{(2)}\supset \cdots \supset \frak{h}^{(N)}=\left\{
0\right\}
\end{equation*}
for some $N\geq 2$, where $\frak{h}^{(1)}=\frak{h}$ , $\frak{h}^{(k)}=\left[ 
\frak{h}^{(k-1)},\frak{h}^{(k-1)}\right] $ for $k>1$, and each containment $%
\supset $ is proper. From Corollary \ref{coro-ideals}, $J=\frak{h}\times 
\frak{h}\times \frak{h}$ is an ideal of $A(\frak{g})$. Now 
\begin{equation*}
J^{2}=\frak{h}^{(2)}\times \frak{h}^{(2)}\times \frak{h}^{(2)}\subset \frak{h%
}\times \frak{h}\times \frak{h}=J,
\end{equation*}
is a proper containment, contradicting (\ref{Jsqrd=J}). Thus $\frak{g}$ has
no nonzero solvable ideals, which implies that $\frak{g}$ is semisimple.
\end{proof}

\section{The Radical; Levi-Malcev Decompositions}

The \emph{radical} of an algebra $\mathcal{A}$, denoted by $\mathrm{rad}\ 
\mathcal{A}$,$\mathrm{\ }$is an ideal of $\mathcal{A}$ characterized as
follows: if $\mathcal{J}$ is an ideal of $\mathcal{A}$ and $\mathcal{A}/%
\mathcal{J}$ is semisimple, then $\mathrm{rad}\ \mathcal{A}\subseteq 
\mathcal{J}$. The existence of $\mathrm{rad}\ \mathcal{A}$ can be shown
using the Chinese Remainder Theorem; see Walcher \cite{wbook}, pp. 3-7. The
radical of a Lie algebra can also be defined as being its maximal solvable
ideal of $\frak{g}$ \cite{jacobson} \cite{sw}. For an algebra $\mathcal{A}$,
if there exists a semisimple subalgebra $\mathcal{S}$ such that $\mathcal{A}=%
\mathcal{S}\oplus \mathrm{rad}\ \mathcal{A}$, then $\mathcal{A}$ is said to
have a \emph{Levi-Malcev decomposition} and $\mathcal{S}$ is called a \emph{%
Levi factor}. For instance, every Lie algebra has a Levi-Malcev
decomposition \cite{jacobson}. For a commutative algebra with a Levi-Malcev
decomposition, the associated quadratic differential equation can be
decoupled into an equation in the Levi factor and a nonautonomous equation
in the radical \cite{ks1}.

The relationship between the radical of a Lie algebra and the radical of its
Nahm algebra is contained in the following result.

\begin{theorem}
\label{thm-radical}$\mathrm{rad}\ A(\frak{g})=A(\mathrm{rad}\ \frak{g}).$
\end{theorem}

\begin{proof}
First observe that 
\begin{equation}
A(\frak{g})/A(\mathrm{rad}\ \frak{g})\cong A(\frak{g}/\mathrm{rad}\ \frak{g})
\label{A(g)/A(r)=A(g/r)}
\end{equation}
by a standard isomorphism theorem for vector spaces. Thus $A(\frak{g})/A(%
\mathrm{rad}\ \frak{g})$ is semisimple by Theorem \ref{thm-semisimple}. Now
suppose $J\subseteq A(\mathrm{rad}\ \frak{g})$ is an ideal of $\frak{g}$
such that $A(\frak{g})/J$ is semisimple. Then $A(\mathrm{rad}\ \frak{g})/J$
is an ideal of $A(\frak{g})/J$. Now $A(\mathrm{rad}\ \frak{g})$ is a
solvable ideal of $A(\frak{g})$ since $\mathrm{rad}\ \frak{g}$ is a solvable
ideal of $\frak{g}$. This implies $A(\mathrm{rad}\ \frak{g})/J$ is a
solvable ideal of $A(\frak{g})/J$. Since $A(\frak{g})/J$ is semisimple, we
must have $J=A(\mathrm{rad}\ \frak{g})$. This completes the proof.
\end{proof}

Since every Lie algebra $\frak{g}$ has a Levi-Malcev decomposition $\frak{g}=%
\frak{s}\oplus \mathrm{rad}\ \frak{g}$ where $\frak{s}$ is semisimple, we
immediately obtain the following result for Nahm algebras.

\begin{corollary}
\label{coro-levi}Every Nahm algebra has a Levi-Malcev decomposition.
\end{corollary}

\begin{proof}
Let $\frak{g}=\frak{s}\oplus \mathrm{rad}\ \frak{g}$ be a Levi-Malcev
decomposition with $\frak{s}$ a semisimple Levi factor. Then $A(\frak{g})=A(%
\frak{s})\oplus A(\mathrm{rad}\ \frak{g})$. By Theorem \ref{thm-semisimple}, 
$A(\frak{s})$ is semisimple. By Theorem \ref{thm-radical}, $A(\mathrm{rad}\ 
\frak{g})=\mathrm{rad}\ A(\frak{g})$. This completes the proof.
\end{proof}

\section{Invariant Bilinear Forms}

For many interesting classes of algebras $\mathcal{A}$, there exists a
bilinear form $F:\mathcal{A}\times \mathcal{A}\rightarrow \mathbb{K}$ which
reflects the structure of the algebra. The forms of interest are those which
are \emph{invariant} (or \emph{associative}), meaning 
\begin{equation*}
F(X\cdot Y,Z)=F(X,Y\cdot Z)
\end{equation*}
for all $X,Y,Z\in \mathcal{A}$. The \emph{radical} of a bilinear form $F:%
\mathcal{A}\times \mathcal{A}\rightarrow \mathbb{K}$ on an algebra $\mathcal{%
A}$ is the subspace $\mathrm{rad\ }F=\left\{ X\in \mathcal{A}:F(X,Y)=0\text{
for all }Y\in \mathcal{A}\right\} $. If $F$ is invariant, then the radical
is an ideal: for $Z\in \mathrm{rad\ }F$ and for all $X,Y\in \mathcal{A}$, $%
F(Y,X\cdot Z)=F(Y\cdot X,Z)=0$, which implies $X\cdot Z\in \mathrm{rad\ }F$.
A bilinear form $F:\mathcal{A}\times \mathcal{A}\rightarrow \mathbb{K}$ is 
\emph{nondegenerate} if its radical is the zero subspace, i.e., $\mathrm{%
rad\ }F=\left\{ 0\right\} $.

If a commutative algebra $\mathcal{A}$ has a nondegenerate invariant form $F$%
, then the associated quadratic equation $\dot{X}=X^{2}$ in $\mathcal{A}$
turns out to be a gradient equation with potential function $\phi (X)=\frac{1%
}{3}F(X,X^{2})$, that is, $X^{2}=(\nabla \phi )(X)$ for all $X\in \mathcal{A}
$. Conversely, if the vector field $X\longmapsto X^{2}$ has a potential
function $\phi :\mathcal{A}\rightarrow \mathbb{K}$ satisfying $X^{2}=(\nabla
\phi )(X)$ for all $X\in \mathcal{A}$, then the bilinear form $F:\mathcal{A}%
\times \mathcal{A}\rightarrow \mathbb{K}$ defined by 
\begin{equation*}
F(X,Y)=\frac{1}{2}\phi ^{(2)}(0)(X,Y)
\end{equation*}
is invariant and nondegenerate. (Note that we are using the term
``gradient'' in a broad sense, for we are not requiring that the form $F$ be
positive definite. Even in the nondegenerate case, one still has the
property that trajectories of the differential equation cross the quadrics $%
F(X,X)=0$ orthogonally relative to $F$ itself.) For further details, see
Walcher \cite{wbook}, p.80ff.

Let $\frak{g}$ be a Lie algebra, and let $\rho :\frak{g}\rightarrow gl(V)$
be a representation of $\frak{g}$ as a Lie algebra of linear transformations
on some finite-dimensional vector space $V$. The \emph{trace form} of $\frak{%
g}$ induced by $\rho $ is the bilinear form $B_{\rho }:\frak{g}\times \frak{g%
}\rightarrow \mathbb{K}$ defined by 
\begin{equation}
B_{\rho }(x,y)=\mathrm{tr}(\rho (x)\rho (y))  \label{rho-form}
\end{equation}
for $x,y\in \frak{g}$ where $\mathrm{tr}$ denotes the trace. Using the fact
that $\rho $ is a representation, it is easy to show that $B_{\rho }$ is
invariant. The most important trace form on a Lie algebra is the one induced
by the adjoint representation, which is called the \emph{Killing form} $%
\kappa :\frak{g}\times \frak{g}\rightarrow \mathbb{K}$, and is given by 
\begin{equation}
\kappa (x,y)=\mathrm{tr}((\mathrm{ad}\ x)(\mathrm{ad}\ y))  \label{killing}
\end{equation}
for $x,y\in \frak{g}$. The Killing form measures the structure of $\frak{g}$
in the following sense: $\frak{g}$ is semisimple if and only if $\kappa $ is
nondegenerate \cite{jacobson} \cite{sw}.

We now introduce a related bilinear form on the Nahm algebra $A(\frak{g})$
which will turn out to measure the structure of $A(\frak{g})$ in a similar
way.

\begin{definition}
\label{defn-trace-form}Let $\frak{g}$ be a Lie algebra, and let $\rho :\frak{%
g}\rightarrow gl(V)$ be a representation of $\frak{g}$. The \emph{induced
trace form} $C_{\rho }:A\left( \frak{g}\right) \times A\left( \frak{g}%
\right) \rightarrow \mathbb{K}$ on the associated Nahm algebra $A(\frak{g})$
is defined by 
\begin{equation}
C_{\rho }(X,Y)=\mathrm{tr}(L_{\rho }(X)L_{\rho }(Y))  \label{trace-form}
\end{equation}
for $X,Y\in A$. In case $\rho =\mathrm{ad}$, the adjoint representation,
then we simply write $C\equiv C_{\mathrm{ad}}$, and we refer to $C$ as the 
\emph{standard form} on $A(\frak{g})$.
\end{definition}

(See (\ref{L-rho}) for the definition of the operator $L_{\rho }(X)$.)

The induced trace form has an immediate characterization in terms of the
trace form of $\frak{g}$.

\begin{theorem}
\label{thm-trace}For $X,Y\in A$, 
\begin{equation}
C_{\rho }(X,Y)=-\frac{1}{2}\left( B_{\rho }(x_{1},y_{1})+B_{\rho
}(x_{2},y_{2})+B_{\rho }(x_{3},y_{3})\right) \text{.}  \label{C-B}
\end{equation}
\end{theorem}

\begin{proof}
Compute $L_{\rho }(X)L_{\rho }(Y)$ using (\ref{L-rho}) and take the trace.
\end{proof}

\begin{remark}
Recall that in addition to having the structure of a Nahm algebra, $\frak{g}%
\times \frak{g}\times \frak{g}$ is also a Lie algebra with the
component-wise bracket. Since any representation $\rho :\frak{g}\rightarrow
gl(V)$ trivially lifts to a representation $\hat{\rho}:\frak{g}\times \frak{g%
}\times \frak{g}\rightarrow gl(V\times V\times V)$, trace forms $B_{\rho }$
on $\frak{g}$ induce trace forms $B_{\hat{\rho}}$ on $\frak{g}\times \frak{g}%
\times \frak{g}$ by 
\begin{equation*}
B_{\hat{\rho}}(X,Y)=B_{\rho }(x_{1},y_{1})+B_{\rho }(x_{2},y_{2})+B_{\rho
}(x_{3},y_{3})
\end{equation*}
for $X,Y\in \frak{g}\times \frak{g}\times \frak{g}$. We thus have the
following relationship between the induced trace form $C_{\rho }$ on the
Nahm algebra $A(\frak{g})$ and the induced trace form $B_{\hat{\rho}}$ on
the Lie algebra $\frak{g}\times \frak{g}\times \frak{g}$: 
\begin{equation*}
C_{\rho }(X,Y)=-\frac{1}{2}B_{\hat{\rho}}(X,Y)
\end{equation*}
for $X,Y\in A(\frak{g})$.
\end{remark}

\begin{theorem}
\label{thm-invariant-form}The induced trace form $C_{\rho }:A\times
A\rightarrow \mathbb{K}$ is invariant, i.e., $C_{\rho }(XY,Z)=C_{\rho
}(X,YZ) $ for all $X,Y,Z\in A$.
\end{theorem}

\begin{proof}
Let $X,Y,Z\in A$ be given. Identify indices modulo $3$: $x_{i}=x_{i+3}$,
etc. Using (\ref{C-B}) and the invariance of the trace form $B_{\rho }:\frak{%
g}\times \frak{g}\rightarrow \mathbb{K}$, we compute 
\begin{eqnarray*}
-2C_{\rho }(XY,Z) &=&\sum_{i=1}^{3}B_{\rho }\left( \frac{1}{2}\left(
[x_{i+1},y_{i+2}]+[x_{i+2},y_{i+1}]\right) ,z_{i}\right) \\
&=&\sum_{i=1}^{3}B_{\rho }\left( x_{i+1},\frac{1}{2}[y_{i+2},z_{i}]\right)
+\sum_{i=1}^{3}B_{\rho }\left( x_{i+2},\frac{1}{2}[y_{i+1},z_{i}]\right) .
\end{eqnarray*}
Since we are summing over all terms, we may reindex and combine them to
obtain 
\begin{eqnarray*}
-2C_{\rho }(XY,Z) &=&\sum_{i=1}^{3}B_{\rho }\left( x_{i},\frac{1}{2}[%
y_{i+1},z_{i+2}]\right) +\sum_{i=1}^{3}B_{\rho }\left( x_{i},\frac{1}{2}[%
y_{i+2},z_{i+1}]\right) \\
&=&\sum_{i=1}^{3}B_{\rho }\left( x_{i},\frac{1}{2}\left(
[y_{i+1},z_{i+2}]+[y_{i+2},z_{i+1}]\right) \right) \\
&=&-2C_{\rho }(X,YZ).
\end{eqnarray*}
This completes the proof.
\end{proof}

For a Nahm algebra $A(\frak{g})$, the radical of an induced trace form is,
of course, related to the radical of the corresponding trace form on $\frak{g%
}$.

\begin{theorem}
\label{thm-rad}$\mathrm{rad\ }C_{\rho }=\mathrm{rad\ }B_{\rho }\times 
\mathrm{rad\ }B_{\rho }\times \mathrm{rad\ }B_{\rho }$.
\end{theorem}

\begin{proof}
Fix $Y=(y_{1},y_{2},y_{3})^{T}\in \mathrm{rad\ }C_{\rho }$. For all $x\in 
\frak{g}$, 
\begin{equation*}
B_{\rho }(x,y_{1})=C_{\rho }\left( \left[ 
\begin{array}{c}
x \\ 
0 \\ 
0
\end{array}
\right] ,\left[ 
\begin{array}{c}
y_{1} \\ 
y_{2} \\ 
y_{3}
\end{array}
\right] \right) =0\text{,}
\end{equation*}
which shows that $y_{1}\in \mathrm{rad\ }B_{\rho }$, and similar
computations show $y_{2},y_{3}\in \mathrm{rad\ }B_{\rho }$. The reverse
inclusion is clear.
\end{proof}

The following corollaries are immediate.

\begin{corollary}
\label{coro-nondegen}An induced form $C_{\rho }:A\times A\rightarrow \mathbb{%
K}$ \ is nondegenerate if and only if the trace form $B_{\rho }:\frak{g}%
\times \frak{g}\rightarrow \mathbb{K}$ \ is nondegenerate.
\end{corollary}

\begin{corollary}
\label{coro-nondeg-semi}The following are equivalent.

\begin{enumerate}
\item  The standard form $C:A\times A\rightarrow \mathbb{K}$ \ is
nondegenerate.

\item  The Killing form $\kappa :\frak{g}\times \frak{g}\rightarrow \mathbb{K%
}$ \ is nondegenerate.

\item  $\frak{g}$ is semisimple.

\item  $A\left( \frak{g}\right) $ is semisimple.
\end{enumerate}
\end{corollary}

Recall that a Lie algebra $\frak{g}$ is said to be \emph{compact} if its
associated (connected) Lie group is compact. A semisimple Lie algebra is
compact if and only if its Killing form is negative definite.

\begin{definition}
\label{defn-compact}A Nahm algebra $A\left( \frak{g}\right) $ is said to be 
\emph{compact} if its underlying Lie algebra $\frak{g}$ is compact.
\end{definition}

\begin{theorem}
\label{thm-semi-compact}A semisimple Nahm algebra is compact if and only if
its standard form is positive definite.
\end{theorem}

\begin{proof}
This is immediate from (\ref{C-B}).
\end{proof}

Recalling our earlier discussion of gradients, it follows that the Nahm
equations (\ref{eqn1})-(\ref{eqn3}) in a compact semisimple Nahm algebra
form a gradient system in the traditional sense.

\begin{remark}
Recall the diagonal subalgebra $\Delta \left( \frak{g}\right) $ of a Nahm
algebra $A\left( \frak{g}\right) $. For $\Delta (x)\in \Delta \left( \frak{g}%
\right) $, $Y\in A\left( \frak{g}\right) $, we have from (\ref{C-B}) 
\begin{equation*}
C(\Delta (x),Y)=-\frac{1}{2}B_{\rho }(x,y_{1}+y_{2}+y_{3})\text{.}
\end{equation*}
Thus the orthogonal complement of $\Delta \left( \frak{g}\right) $ relative
to the standard form is the subspace 
\begin{equation*}
W_{\mathrm{rad}}\left( \frak{g}\right) =\left\{ Y\in A\left( \frak{g}\right)
:y_{1}+y_{2}+y_{3}\in \mathrm{rad\ }\kappa \right\} \text{.}
\end{equation*}
The intersection of this subspace with the diagonal subalgebra is 
\begin{equation*}
\Delta \left( \frak{g}\right) \cap W_{\mathrm{rad}}\left( \frak{g}\right)
=\left\{ \Delta (x)\in \Delta \left( \frak{g}\right) :x\in \mathrm{rad\ }%
\kappa \right\} \text{.}
\end{equation*}
This is simply a copy of $\mathrm{rad\ }\kappa $ itself. In particular, we
see from Definition \ref{defn-diagonal} that 
\begin{equation*}
W_{\mathrm{rad}}\left( \frak{g}\right) =W\left( \frak{g}\right) \text{ if
and only if }\mathrm{rad\ }\kappa =\left\{ 0\right\} \text{,}
\end{equation*}
that is, if and only if $\kappa $ and $C$ are nondegenerate (Corollary \ref
{coro-nondeg-semi}).
\end{remark}

\section{Derivations}

A \emph{derivation} of an algebra $\mathcal{A}$ is a linear transformation $%
D:\mathcal{A}\rightarrow \mathcal{A}$ satisfying $D\left( XY\right) =\left(
DX\right) Y+X\left( DY\right) $ for all $X,Y\in \mathcal{A}$. Let $\mathrm{%
Der}\left( \mathcal{A}\right) $ denote the space of all derivations of $%
\mathcal{A}$; this is a Lie subalgebra of $gl\left( \mathcal{A}\right) $.

If $\mathcal{A}$ is commutative, then derivations of $\mathcal{A}$ are
linear infinitesimal symmetries of the quadratic differential equation $\dot{%
X}=X^{2}$ in $\mathcal{A}$. For $D\in \mathrm{Der}\left( \mathcal{A}\right) $%
, $DX(t;P)=\triangledown X(t;P)\cdot DP$, where $\triangledown $ represents
the derivative with respect to the $\mathcal{A}$-variables. If $D\in \mathrm{%
Der}\left( \mathcal{A}\right) $ and $P\in \mathcal{A}$ are such that $%
DP=P^{2}$, then $X\left( t\right) =e^{tD}P$ turns out to be the unique
solution with initial value $P$. For more details, see Walcher \cite{wbook},
Kinyon and Sagle \cite{ks1} \cite{ks2} \cite{ks3}.

Now for $1\leq i<j\leq 3$, let $E_{ij}=\mathbf{e}_{i}\mathbf{e}_{j}^{t}-%
\mathbf{e}_{j}\mathbf{e}_{i}^{t}$. Then $\left\{ E_{ij}\right\} _{1\leq
i<j\leq 3}$ is a basis for the Lie algebra $so\left( 3,\mathbb{K}\right) $
of $3\times 3$ skew-symmetric matrices.

\begin{theorem}
\label{thm-so(3)-der}$so\left( 3,\mathbb{K}\right) $ is a Lie subalgebra of $%
\mathrm{Der}\left( A\left( \frak{g}\right) \right) $.
\end{theorem}

\begin{proof}
We will show that $E_{12}$ is a derivation. That $E_{13}$ and $E_{23}$ are
derivations follow similarly. For $X=\left( x_{1},x_{2},x_{3}\right)
^{T},Y=\left( y_{1},y_{2},y_{3}\right) ^{T}\in A\left( \frak{g}\right) $, we
compute 
\begin{eqnarray}
E_{12}(XY) &=&\left( 
\begin{array}{ccc}
0 & 1 & 0 \\ 
-1 & 0 & 0 \\ 
0 & 0 & 0
\end{array}
\right) \left( \frac{1}{2}\left( 
\begin{array}{c}
\lbrack x_{2},y_{3}]+[y_{2},x_{3}] \\ 
\lbrack x_{3},y_{1}]+[y_{3},x_{1}] \\ 
\lbrack x_{1},y_{2}]+[y_{1},x_{2}]
\end{array}
\right) \right)  \notag \\
&=&\frac{1}{2}\left( 
\begin{array}{c}
\lbrack x_{3},y_{1}]+[y_{3},x_{1}] \\ 
-[x_{2},y_{3}]-[y_{2},x_{3}] \\ 
0
\end{array}
\right) \text{.}  \label{der1}
\end{eqnarray}
On the other hand, we have 
\begin{eqnarray}
\left( E_{12}X\right) Y+X\left( E_{12}Y\right) &=&\left( 
\begin{array}{c}
x_{2} \\ 
-x_{1} \\ 
0
\end{array}
\right) \left( 
\begin{array}{c}
y_{1} \\ 
y_{2} \\ 
y_{3}
\end{array}
\right) +\left( 
\begin{array}{c}
x_{1} \\ 
x_{2} \\ 
x_{3}
\end{array}
\right) \left( 
\begin{array}{c}
y_{2} \\ 
-y_{1} \\ 
0
\end{array}
\right)  \notag \\
&=&\frac{1}{2}\left( 
\begin{array}{c}
-\left[ x_{1},y_{3}\right] -\left[ y_{1},x_{3}\right] \\ 
\lbrack y_{3},x_{2}]+[x_{3},y_{2}] \\ 
\left[ x_{2},y_{2}\right] -\left[ y_{1},x_{1}\right] -\left[ x_{1},y_{1}%
\right] +\left[ y_{2},x_{2}\right]
\end{array}
\right) .  \label{der2}
\end{eqnarray}
Comparing (\ref{der1}) and (\ref{der2}), and using the skew-symmetry of the
Lie bracket, the result follows.
\end{proof}

Next we identify another subalgebra of $\mathrm{Der}(A\left( \frak{g}\right)
)$. Observe that the mapping $\mathrm{diag}(\mathrm{ad\ }(\cdot )):\frak{g}%
\rightarrow gl(\mathcal{A})$ is an isomorphic copy of the adjoint
representation of $\frak{g}$. Let 
\begin{equation}
\mathrm{diag}(\mathrm{ad}(\frak{g))}=\left\{ \mathrm{diag}(\mathrm{ad\ }%
x):x\in \frak{g}\right\} .  \label{ad(D(g))}
\end{equation}
This is an isomorphic copy of $\mathrm{ad}\left( \frak{g}\right) $, and, if
the adjoint representation is faithful, of $\frak{g}$ itself.

\begin{theorem}
\label{thm-adD(g)-der}$\mathrm{diag}(\mathrm{ad}(\frak{g))}$ is a Lie
subalgebra of $\mathrm{Der}(A\left( \frak{g}\right) )$.
\end{theorem}

\begin{proof}
The identity 
\begin{equation*}
\mathrm{diag}(\mathrm{ad\ }x)(YZ)=(\mathrm{diag}(\mathrm{ad\ }x)Y)Z+Y(%
\mathrm{diag}(\mathrm{ad\ }x)Z)
\end{equation*}
is an easy consequence of the Jacobi identity in $\frak{g}$.
\end{proof}

Obviously $\mathrm{diag}(\mathrm{ad}(\frak{g))}\cap so\left( 3,\mathbb{K}%
\right) =\left\{ 0\right\} $, and thus $\mathrm{Der}(A\left( \frak{g}\right)
)$ contains $\mathrm{diag}(\mathrm{ad}(\frak{g))}\oplus so\left( 3,\mathbb{K}%
\right) $ as a direct sum of vector spaces. In addition, for $x\in \frak{g}$
and $M\in so\left( 3,\mathbb{K}\right) $, a direct calculation yields 
\begin{equation}
\left[ \mathrm{diag}(\mathrm{ad\ }x),M\right] =0.  \label{commute}
\end{equation}
Thus $\mathrm{diag}(\mathrm{ad}(\frak{g))}\oplus so\left( 3,\mathbb{K}%
\right) $ is also an internal direct sum of Lie subalgebras of $\mathrm{Der}%
(A\left( \frak{g}\right) )$.

The rest of this section will be devoted to proving the following result.

\begin{theorem}
\label{thm-simple-der}Let $A\left( \frak{g}\right) $ be a simple Nahm
algebra. Then 
\begin{equation*}
\mathrm{Der}(A\left( \frak{g}\right) )=\mathrm{diag}(\mathrm{ad}(\frak{g))}%
\oplus so\left( 3,\mathbb{K}\right) .
\end{equation*}
\end{theorem}

Note that the result will turn out to hold for both $\mathbb{K}=\mathbb{C}$
and $\mathbb{K}=\mathbb{R}$.

In the course of the discussion that follows, we will have frequent occasion
to use the equivalence of the simplicity of $A\left( \frak{g}\right) $ with
the simplicity of $\frak{g}$ (Theorem \ref{thm-simple}) without explicitly
mentioning it.

Recall the left multiplication operator $L\left( X\right) :A\left( \frak{g}%
\right) \rightarrow A\left( \frak{g}\right) $ given by (\ref{L(X)}). We will
denote the identity operator in $\frak{g}$ or $A\left( \frak{g}\right) $ by $%
I$, and let the context clarify which is meant. We begin with a version of
Schur's Lemma for complex, simple Nahm algebras.

\begin{lemma}
\label{lemma-schur}Let $A\left( \frak{g}\right) $ be a complex, simple Nahm
algebra, and let $T\in gl(A\left( \frak{g}\right) )$ satisfy $T\circ L\left(
X\right) =L\left( X\right) \circ T$ for all $X\in A\left( \frak{g}\right) $.
Then there exists $\lambda \in \mathbb{C}$ such that $T=\lambda I$.
\end{lemma}

\begin{proof}
For $x\in \frak{g}$, let $X=\left( x,0,0\right) ^{t}$. Multiplying matrices,
we find that the equation $T\circ L\left( X\right) =L\left( X\right) \circ T$
is equivalent to the equations 
\begin{eqnarray*}
-T_{12}\circ \mathrm{ad\ }x &=&T_{13}\circ \mathrm{ad\ }x=0 \\
-\mathrm{ad\ }x\circ T_{31} &=&\mathrm{ad\ }x\circ T_{21}=0 \\
T_{23}\circ \mathrm{ad\ }x &=&-\mathrm{ad\ }x\circ T_{32} \\
-T_{23}\circ \mathrm{ad\ }x &=&-\mathrm{ad\ }x\circ T_{33} \\
T_{33}\circ \mathrm{ad\ }x &=&\mathrm{ad\ }x\circ T_{22} \\
-T_{32}\circ \mathrm{ad\ }x &=&\mathrm{ad\ }x\circ T_{23}.
\end{eqnarray*}
If we similarly let $X=\left( 0,x,0\right) ^{t}$ and $X=\left( 0,0,x\right)
^{t}$, and consider the corresponding matrix equations, then by matching
matrix entries, we finally obtain the following system of equations in $%
\frak{g}$: 
\begin{eqnarray}
T_{ij}\circ \mathrm{ad\ }x &=&\mathrm{ad\ }x\circ T_{ij}=0  \label{schur1} \\
T_{ii}\circ \mathrm{ad\ }x &=&\mathrm{ad\ }x\circ T_{jj}  \label{schur2}
\end{eqnarray}
for $i\neq j$. By Schur's Lemma, (\ref{schur1}) implies $T_{ij}=0$ for $%
i\neq j$. Identifying indices modulo $3$, (\ref{schur2}) implies 
\begin{equation*}
T_{ii}\circ \mathrm{ad\ }x=\mathrm{ad\ }x\circ T_{i+1,i+1}=T_{i+2,i+2}\circ 
\mathrm{ad\ }x=\mathrm{ad\ }x\circ T_{ii}
\end{equation*}
for $i=1,2,3$. By Schur's Lemma, for $i=1,2,3$, there exists $\lambda
_{i}\in \mathbb{C}$ such that $T_{ii}=\lambda _{i}I$. But then (\ref{schur2}%
) implies that $\lambda _{1}=\lambda _{2}=\lambda _{3}$. This completes the
proof.
\end{proof}

For a Nahm algebra $A\left( \frak{g}\right) $, recall the standard form $%
C:A\left( \frak{g}\right) \times A\left( \frak{g}\right) \rightarrow \mathbb{%
K}$ given by (\ref{trace-form}) or, equivalently, (\ref{C-B}) where $\rho $
is the adjoint representation of $\frak{g}$. Assume $\frak{g}$ is semisimple
so that $C$ is nondegenerate (Corollary \ref{coro-nondeg-semi}). For a
linear transformation $T:A\left( \frak{g}\right) \rightarrow A\left( \frak{g}%
\right) $, let $T^{c}:A\left( \frak{g}\right) \rightarrow A\left( \frak{g}%
\right) $ denote the $C$-transpose of $T$ defined by 
\begin{equation}
C(T^{c}X,Y)=C(X,TY)  \label{transpose}
\end{equation}
for all $X,Y\in A\left( \frak{g}\right) $.

\begin{lemma}
\label{lemma-T+Tc}Let $A\left( \frak{g}\right) $ be a complex, simple Nahm
algebra, and let $T\in \mathrm{Der}\left( A\left( \frak{g}\right) \right) $.
Then there exists $\lambda \in \mathbb{C}$ such that $T+T^{c}=\lambda I$.
\end{lemma}

\begin{proof}
For $X,Y,Z\in A\left( \frak{g}\right) $, we compute 
\begin{eqnarray*}
C(X\left( T^{c}Y\right) ,Z) &=&C(T^{c}Y,XZ) \\
&=&C(Y,T(XZ)) \\
&=&C(Y,(TX)Z+X(TZ)) \\
&=&C(Y,(TX)Z)+C(Y,X(TZ)) \\
&=&C(Y(TX),Z)+C(YX,TZ) \\
&=&C(Y(TX),Z)+C(T^{c}(YX),Z) \\
&=&C(Y(TX)+T^{c}(YX),Z)
\end{eqnarray*}
using the invariance and bilinearity of the standard form, (\ref{transpose}%
), and $T\in \mathrm{Der}(A\left( \frak{g}\right) )$. Since $C$ is
nondegenerate (Corollary \ref{coro-nondeg-semi}), 
\begin{equation}
X(T^{c}Y)=Y(TX)+T^{c}(YX)  \label{c-der1}
\end{equation}
for all $X,Y\in A\left( \frak{g}\right) $. On the other hand, since $T$ is a
derivation, 
\begin{equation}
X(TY)=-Y(TX)+T(XY)  \label{c-der2}
\end{equation}
for all $X,Y\in A\left( \frak{g}\right) $. Adding (\ref{c-der1}) and (\ref
{c-der2}), we obtain 
\begin{equation*}
L(X)(T+T^{c})Y=(T+T^{c})L(X)Y
\end{equation*}
for all $X,Y\in A\left( \frak{g}\right) $. By Lemma \ref{lemma-schur}, there
exists $\lambda \in \mathbb{C}$ such that $T+T^{c}=\lambda I$.
\end{proof}

It is easy to see that $\mathrm{diag}(\mathrm{ad\ }x)$ is $C$-skew symmetric
for all $x\in \frak{g}$, and any matrix in $so(3,\mathbb{C})$ is $C$-skew
symmetric as a linear transformation on $A$. Thus once Theorem \ref
{thm-simple-der} is established, it will follow that the conclusion of Lemma 
\ref{lemma-T+Tc} can be strengthened to the assertion that $T$ is $C$-skew
symmetric. The conclusion of the lemma as it is presently stated will be
used in the proof of Theorem \ref{thm-simple-der}.

Let $A\left( \frak{g}\right) $ be a Nahm algebra and let $T\in gl(A\left( 
\frak{g}\right) )$ have the usual block matrix representation $T=\left[
T_{ij}\right] $. Define linear transformations $T_{\mathrm{diag}},T_{\mathrm{%
off}}\in gl(A\left( \frak{g}\right) )$ by 
\begin{equation}
T_{\mathrm{diag}}=\left( 
\begin{array}{ccc}
T_{11} & 0 & 0 \\ 
0 & T_{22} & 0 \\ 
0 & 0 & T_{33}
\end{array}
\right)  \label{T-diag}
\end{equation}
and 
\begin{equation}
T_{\mathrm{off}}=\left( 
\begin{array}{ccc}
0 & T_{12} & T_{13} \\ 
T_{21} & 0 & T_{23} \\ 
T_{31} & T_{32} & 0
\end{array}
\right) .  \label{T-off}
\end{equation}
For $X\in A\left( \frak{g}\right) $, we compute 
\begin{equation}
T(X^{2})=\left( 
\begin{array}{c}
T_{11}\left[ x_{2},x_{3}\right] +T_{12}\left[ x_{3},x_{1}\right] +T_{13}%
\left[ x_{1},x_{2}\right] \\ 
T_{21}\left[ x_{2},x_{3}\right] +T_{22}\left[ x_{3},x_{1}\right] +T_{23}%
\left[ x_{1},x_{2}\right] \\ 
T_{31}\left[ x_{2},x_{3}\right] +T_{32}\left[ x_{3},x_{1}\right] +T_{33}%
\left[ x_{1},x_{2}\right]
\end{array}
\right)  \label{T1}
\end{equation}
and $2\left( TX\right) X=$ 
\begin{equation}
\left( 
\begin{array}{c}
\left[ T_{21}x_{1}+T_{22}x_{2}+T_{23}x_{3},x_{3}\right] +\left[
x_{2},T_{31}x_{1}+T_{32}x_{2}+T_{33}x_{3}\right] \\ 
\left[ T_{31}x_{1}+T_{32}x_{2}+T_{33}x_{3},x_{1}\right] +\left[
x_{3},T_{11}x_{1}+T_{12}x_{2}+T_{13}x_{3}\right] \\ 
\left[ T_{11}x_{1}+T_{12}x_{2}+T_{13}x_{3},x_{2}\right] +\left[
x_{1},T_{21}x_{1}+T_{22}x_{2}+T_{23}x_{3}\right]
\end{array}
\right) .  \label{T2}
\end{equation}

\begin{lemma}
\label{lemma-der-split}Let $A\left( \frak{g}\right) $ be a Nahm algebra and
let $T\in \mathrm{Der}(A\left( \frak{g}\right) )$. Then $T_{\mathrm{diag}%
}\in \mathrm{Der}(A\left( \frak{g}\right) )$ and $T_{\mathrm{off}}\in 
\mathrm{Der}(A\left( \frak{g}\right) )$.
\end{lemma}

\begin{proof}
Since $T$ is a derivation, (\ref{T1}) and (\ref{T2}) are equal. Take $%
x_{2}=x_{3}=0$, $x_{1}=x$ in (\ref{T1}) and (\ref{T2}) and match the
entries. This gives $\left[ T_{31}x,x\right] =0$ and $[x,T_{21}x]=0$ for all 
$x\in \frak{g}$. By similar arguments, we obtain 
\begin{equation}
\lbrack T_{ij}x,x]=0  \label{T3}
\end{equation}
for all $x\in \frak{g}$ where $i\neq j$. Linearizing (\ref{T3}) and
rearranging, we have 
\begin{equation}
\lbrack T_{ij}x,y]=[x,T_{ij}y]  \label{lin}
\end{equation}
for all $x,y\in \frak{g}$ where $i\neq j$. Take $x_{1}=0$ in (\ref{T1}) and (%
\ref{T2}), simplify using (\ref{T3}), and match the first entries. This
gives 
\begin{equation}
T_{11}\left[ x_{2},x_{3}\right] =\left[ T_{22}x_{2},x_{3}\right] +\left[
x_{2},T_{33}x_{3}\right]  \label{diagder1}
\end{equation}
for all $x_{2},x_{3}\in \frak{g}$. Successively taking $x_{2}=0$ and $%
x_{3}=0 $ give the equations 
\begin{eqnarray}
T_{22}\left[ x_{3},x_{1}\right] &=&\left[ T_{33}x_{3},x_{1}\right] +\left[
x_{3},T_{11}x_{1}\right]  \label{diagder2} \\
T_{33}\left[ x_{1},x_{2}\right] &=&\left[ T_{11}x_{1},x_{2}\right]
+[x_{1},T_{22}x_{2}]  \label{diagder3}
\end{eqnarray}
for all $x_{1},x_{2},x_{3}\in \frak{g}$. Taken together, (\ref{diagder1}), (%
\ref{diagder2}) and (\ref{diagder3}) imply 
\begin{equation*}
T_{\mathrm{diag}}\left( X^{2}\right) =2X\left( T_{\mathrm{diag}}X\right)
\end{equation*}
for all $X\in A\left( \frak{g}\right) $, i.e., $T_{\mathrm{diag}}$ is a
derivation of $A\left( \frak{g}\right) $. Since $T_{\mathrm{off}}=T-T_{%
\mathrm{diag}}$, $T_{\mathrm{off}}$ is also a derivation.
\end{proof}

\begin{lemma}
\label{lemma-der-off-case}Let $A\left( \frak{g}\right) $ be a complex,
simple Nahm algebra, and let $T\in \mathrm{Der}(A\left( \frak{g}\right) )$
be given. Assume $T=T_{\mathrm{off}}$. Then the action of $T$ on $A\left( 
\frak{g}\right) $ is given by the action of a matrix in $so(3,\mathbb{K})$.
\end{lemma}

\begin{proof}
Using (\ref{T3}), the equality of (\ref{T1}) and (\ref{T2}) simplifies to 
\begin{equation}
\left( 
\begin{array}{c}
T_{12}\left[ x_{3},x_{1}\right] +T_{13}\left[ x_{1},x_{2}\right] \\ 
T_{21}\left[ x_{2},x_{3}\right] +T_{23}\left[ x_{1},x_{2}\right] \\ 
T_{31}\left[ x_{2},x_{3}\right] +T_{32}\left[ x_{3},x_{1}\right]
\end{array}
\right) =\left( 
\begin{array}{c}
\left[ T_{21}x_{1},x_{3}\right] +\left[ x_{2},T_{31}x_{1}\right] \\ 
\left[ T_{32}x_{2},x_{1}\right] +\left[ x_{3},T_{12}x_{2}\right] \\ 
\left[ T_{13}x_{3},x_{2}\right] +\left[ x_{1},T_{23}x_{3}\right]
\end{array}
\right)  \label{T4}
\end{equation}
for all $x_{1},x_{2},x_{3}\in \frak{g}$. Set $x_{3}=0$ in (\ref{T4}) and
match entries. Using (\ref{lin}), this gives 
\begin{eqnarray*}
T_{13}\left[ x_{1},x_{2}\right] &=&\left[ x_{2},T_{31}x_{1}\right]
=-[x_{1},T_{31}x_{2}] \\
T_{23}\left[ x_{1},x_{2}\right] &=&\left[ T_{32}x_{2},x_{1}\right]
=-[x_{1},T_{32}x_{2}]
\end{eqnarray*}
for all $x_{1},x_{2}\in \frak{g}$. Similar calculations give 
\begin{equation}
T_{ij}[x,y]=-[x,T_{ji}y]  \label{com1}
\end{equation}
for all $x,y\in \frak{g}$ where $i\neq j$. Iterating (\ref{com1}), we have 
\begin{eqnarray*}
T_{ij}[x,[y,z]] &=&-[x,T_{ji}[y,z]] \\
&=&[x,[y,T_{ij}z]]
\end{eqnarray*}
for all $x,y,z\in \frak{g}$. Thus 
\begin{equation}
T_{ij}\circ \mathrm{ad\ }x\circ \mathrm{ad\ }y=\mathrm{ad\ }x\circ \mathrm{%
ad\ }y\circ T_{ij}  \label{T-ad-ad}
\end{equation}
for all $x,y\in \frak{g}$. Reversing the roles of $x$ and $y$ in (\ref
{T-ad-ad}) and subtracting the resulting equation from (\ref{T-ad-ad}), we
obtain 
\begin{equation}
T_{ij}\circ \mathrm{ad\ }[x,y]=\mathrm{ad\ }[x,y]\circ T_{ij}
\label{T-ad[x,y]}
\end{equation}
for all $x,y\in \frak{g}$ since $\mathrm{ad}:\frak{g}\rightarrow \frak{g}$
is a representation. Since $\frak{g}$ is simple, we have $\left[ \frak{g},%
\frak{g}\right] =\frak{g}$ and thus (\ref{T-ad[x,y]}) implies 
\begin{equation*}
T_{ij}\circ \mathrm{ad\ }x=\mathrm{ad\ }x\circ T_{ij}
\end{equation*}
for all $x\in \frak{g}$. By Schur's Lemma, there exists $\lambda _{ij}\in 
\mathbb{C}$ ($i\neq j$) such that $T_{ij}=\lambda _{ij}I$. If we set $%
\lambda _{ii}=0$ for $i=1,2,3$, then the action of $T$ on $A$ is given by
the action of the matrix $\Lambda =[\lambda _{ij}]$. What remains is to show
that $\Lambda \in so(3,\mathbb{C})$, and for this purpose we will use Lemma 
\ref{lemma-T+Tc}. Let $X,Y\in A\left( \frak{g}\right) $ be given. Using (\ref
{transpose}) and (\ref{C-B}), we compute 
\begin{eqnarray*}
C(T^{c}X,Y) &=&C(X,TY) \\
&=&-\frac{1}{2}\sum_{i=1}^{3}B\left( x_{i},\sum_{j=1}^{3}\lambda
_{ij}y_{j}\right) \\
&=&-\frac{1}{2}\sum_{j=1}^{3}B\left( \sum_{i=1}^{3}\lambda
_{ij}x_{i},y_{j}\right) \\
&=&C(\Lambda ^{t}X,Y)
\end{eqnarray*}
where $\Lambda ^{t}$ is the usual transpose of the matrix $\Lambda $. Since $%
C$ is nondegenerate, the action of $T^{c}$ on $A\left( \frak{g}\right) $ is
given by the action of the matrix $\Lambda ^{t}$. By Lemma \ref{lemma-T+Tc},
we have 
\begin{equation*}
\Lambda +\Lambda ^{t}=\mu I
\end{equation*}
for some $\mu \in \mathbb{C}$. But the diagonal entries of $\Lambda $, and
hence $\Lambda ^{t}$, are all zero, and thus $\mu =0$. It follows that $%
\Lambda \in so(3,\mathbb{C})$. This completes the proof.
\end{proof}

\begin{lemma}
\label{lemma-der-diag-case}Let $A\left( \frak{g}\right) $ be a complex,
simple Nahm algebra, and let $T\in \mathrm{Der}(A\left( \frak{g}\right) )$
be given. Assume $T=T_{\mathrm{diag}}$. Then there exists $x\in \frak{g}$
such that $T=\mathrm{diag}(\mathrm{ad\ }x)$.
\end{lemma}

\begin{proof}
The equality of (\ref{T1}) and (\ref{T2}) gives the following system of
equations 
\begin{eqnarray}
T_{11}\left[ x_{2},x_{3}\right] &=&\left[ T_{22}x_{2},x_{3}\right] +\left[
x_{2},T_{33}x_{3}\right]  \label{T5-1} \\
T_{22}\left[ x_{3},x_{1}\right] &=&\left[ T_{33}x_{3},x_{1}\right] +\left[
x_{3},T_{11}x_{1}\right]  \label{T5-2} \\
T_{33}\left[ x_{1},x_{2}\right] &=&\left[ T_{11}x_{1},x_{2}\right] +\left[
x_{1},T_{22}x_{2}\right]  \label{T5-3}
\end{eqnarray}
for all $x_{1},x_{2},x_{3}\in \frak{g}$. Set $x=x_{1}=x_{2}$ and $y=x_{3}$,
and add (\ref{T5-1}) and (\ref{T5-2}) to obtain 
\begin{eqnarray}
(T_{11}-T_{22})[x,y] &=&\left[ T_{22}x,y\right] +\left[ x,T_{33}y\right] +%
\left[ T_{33}y,x\right] +\left[ y,T_{11}x\right]  \notag \\
&=&[y,(T_{11}-T_{22})x]  \label{anticom}
\end{eqnarray}
for all $x,y\in \frak{g}$. Iterating (\ref{anticom}), we obtain 
\begin{eqnarray*}
(T_{11}-T_{22})[z,[y,x]] &=&(T_{11}-T_{22})[[x,y],z] \\
&=&[z,(T_{11}-T_{22})[x,y]] \\
&=&[z,[y,(T_{11}-T_{22})x]]
\end{eqnarray*}
for all $x,y,z\in \frak{g}$. Thus 
\begin{equation}
(T_{11}-T_{22})\circ \mathrm{ad\ }z\circ \mathrm{ad\ }y=\mathrm{ad\ }z\circ 
\mathrm{ad\ }y\circ (T_{11}-T_{22})  \label{T11-T22}
\end{equation}
for all $y,z\in \frak{g}$. Exchanging $y$ and $z$ in (\ref{T11-T22}) and
subtracting the resulting equation from (\ref{T11-T22}) gives 
\begin{equation*}
(T_{11}-T_{22})\circ \mathrm{ad\ }[z,y]=\mathrm{ad\ }[z,y]\circ
(T_{11}-T_{22})
\end{equation*}
for all $y,z\in \frak{g}$ since $\mathrm{ad}:\frak{g}\rightarrow \frak{g}$
is a representation. Now $[\frak{g},\frak{g}]=\frak{g}$ because $\frak{g}$
is simple, and thus $(T_{11}-T_{22})\circ \mathrm{ad\ }z=\mathrm{ad\ }z\circ
(T_{11}-T_{22})$ for all $z\in \frak{g}$. By Schur's Lemma, there exists $%
\lambda _{12}\in \mathbb{C}$ such that $T_{11}-T_{22}=\lambda _{12}I$.
Applying this to (\ref{anticom}), we obtain 
\begin{equation*}
\lambda _{12}[x,y]=[y,\lambda _{12}x]=\lambda _{12}[y,x]=-\lambda _{12}[x,y]
\end{equation*}
for all $x,y\in \frak{g}$. Thus $\lambda _{12}=0$, and hence $T_{11}=T_{22}$%
. Similar arguments, \textit{mutatis mutandis}, show that $T_{22}=T_{33}$.
Now (\ref{T5-1}), say, shows that $T_{ii}$ is a derivation of $\frak{g}$.
Since $\frak{g}$ is simple, there exists $x\in \frak{g}$ such that $T_{ii}=%
\mathrm{ad\ }x$. It follows that $T=\mathrm{diag}(\mathrm{ad\ }x)$ as
claimed.
\end{proof}

Finally, we complete the proof of Theorem \ref{thm-simple-der}. If $\mathbb{K%
}=\mathbb{C}$, then the result follows from Lemmas \ref{lemma-der-split}, 
\ref{lemma-der-off-case}, and \ref{lemma-der-diag-case}. Now suppose $%
\mathbb{K}=\mathbb{R}$. By Theorems \ref{thm-adD(g)-der} and \ref
{thm-so(3)-der}, we have that $\mathrm{diag}(\mathrm{ad}(\frak{g}))\oplus
so(3,\mathbb{R})$ is a subalgebra of $\mathrm{Der}\left( A(\frak{g})\right) $%
. Now the complexification of $A(\frak{g})$ is simple, and its derivation
algebra is the complexification of $\mathrm{Der}\left( A(\frak{g})\right) $.
But the real dimension of a derivation algebra is equal to the complex
dimension of its complexification. This gives us the desired result.

\section{Automorphisms}

An \emph{automorphism} of an algebra $\mathcal{A}$ is an invertible linear
mapping $\phi :\mathcal{A}\rightarrow \mathcal{A}$ satisfying $\phi
(XY)=\phi (X)\phi (Y)$ for all $X,Y\in \mathcal{A}$. Let $\mathrm{Aut}\left( 
\mathcal{A}\right) $ denote the set of all automorphisms of $\mathcal{A}$;
this is a closed (Lie) subgroup of $GL\left( \mathcal{A}\right) $, and $%
\mathrm{Der}\left( \mathcal{A}\right) $ is the Lie algebra of $\mathrm{Aut}%
\left( \mathcal{A}\right) $ (see Sagle and Walde \cite{sw}).

If $\mathcal{A}$ is commutative, then the automorphisms of $\mathcal{A}$ are
the linear symmetries of the quadratic differential equation $\dot{X}=X^{2}$
occurring in $\mathcal{A}$, that is, they are solution preserving. Let $%
X\left( t\right) =X(t;P)$ denote the unique solution with initial value $%
P\in \mathcal{A}$, and let $Y\left( t\right) =\phi \left( X\left( t\right)
\right) $ for $\phi \in \mathrm{Aut}\left( \mathcal{A}\right) $. Then 
\begin{equation*}
\dot{Y}\left( t\right) =\phi \left( \dot{X}\left( t\right) \right) =\phi
\left( X^{2}(t)\right) =\left[ \phi \left( X\left( t\right) \right) \right]
^{2}=Y\left( t\right) ^{2}\text{.}
\end{equation*}
Since $Y\left( 0\right) =\phi P$, we have that $Y\left( t\right) =Y\left(
t;\phi P\right) $ is the unique solution with initial value $\phi P$. For
more on using automorphisms to study quadratic differential equations, see
Walcher \cite{wbook}, Kinyon and Sagle \cite{ks1} \cite{ks2} \cite{ks3},
Hopkins and Kinyon \cite{hk1} \cite{hk2}.

Turning to Nahm algebras, we have the following immediate corollary of
Theorem \ref{thm-so(3)-der}

\begin{corollary}
\label{coro-SO(3)-aut}$SO\left( 3,\mathbb{K}\right) \leq \mathrm{Aut}\left(
A\left( \frak{g}\right) \right) $.
\end{corollary}

\begin{remark}
If $\mathcal{A}$ is a commutative algebra with $\mathbb{Z}_{2}$-grading $%
\mathcal{A}=\mathcal{A}_{0}\oplus \mathcal{A}_{1}$, then the mapping given
by $X_{0}+X_{1}\mapsto X_{0}-X_{1}$ is an automorphism of order $2$, and
conversely, any such automorphism induces a $\mathbb{Z}_{2}$-grading by
setting $\mathcal{A}_{0}$ and $\mathcal{A}_{1}$ equal to the $+1$- and $-1$%
-eigenspaces, respectively. Recall the $\mathbb{Z}_{2}$-grading $A\left( 
\frak{g}\right) =\Delta \left( \frak{g}\right) \oplus W\left( \frak{g}%
\right) $ established in Corollary \ref{coro-grading}. Let $U$ denote the
automorphism that defines the grading. Then $U$ is given by 
\begin{equation*}
U(X)=P_{\Delta }(X)-P_{W}(X)
\end{equation*}
for $X\in A\left( \frak{g}\right) $, where $P_{\Delta }$ and $P_{W}$ are the
projectors defined in (\ref{P_D})\ and (\ref{P_W}), respectively. Computing $%
U(X)$ explicitly in terms of the entries, we find 
\begin{equation*}
U(X)=\frac{1}{3}\left( 
\begin{array}{c}
-x_{1}+2x_{2}+2x_{3} \\ 
2x_{1}-x_{2}+2x_{3} \\ 
2x_{1}+2x_{2}-x_{3}
\end{array}
\right)
\end{equation*}
for all $X\in A\left( \frak{g}\right) $, which implies we can identify $U$
with a $3\times 3$ matrix: 
\begin{equation}
U=\frac{1}{3}\left( 
\begin{array}{ccc}
-1 & 2 & 2 \\ 
2 & -1 & 2 \\ 
2 & 2 & -1
\end{array}
\right) \text{.}  \label{U-auto}
\end{equation}
This matrix is orthogonal and has determinant $1$, and thus $U\in SO(3,%
\mathbb{K})$. A particular derivation $G$ such that $\exp G=U$ is given by 
\begin{equation*}
G=\frac{\pi }{\sqrt{3}}\left( 
\begin{array}{ccc}
0 & 1 & -1 \\ 
-1 & 0 & 1 \\ 
1 & -1 & 0
\end{array}
\right) .
\end{equation*}
\end{remark}

Of course, automorphisms of $\frak{g}$ induce automorphisms of $A(\frak{g})$
in the obvious way: for $\phi \in \mathrm{Aut}(\frak{g})$, we clearly have $%
\mathrm{diag}(\phi )\in \mathrm{Aut}(A(\frak{g}))$.

\begin{proposition}
\label{prop-diagaut}$\mathrm{diag}(\mathrm{Aut}(\frak{g}))\leq \mathrm{Aut}%
(A(\frak{g}))$.
\end{proposition}

Now we show that for simple Nahm algebras, Corollary \ref{coro-SO(3)-aut}
and Proposition \ref{prop-diagaut} describe all the automorphisms. One
particular implication of this is that a Nahm algebra has no outer
automorphisms other than those it inherits from its underlying Lie algebra.

\begin{theorem}
\label{thm-simple-aut}Let $A\left( \frak{g}\right) $ be a simple Nahm
algebra. Then 
\begin{equation*}
\mathrm{Aut}(A\left( \frak{g}\right) )=\mathrm{diag}(\mathrm{Aut}(\frak{g}%
))\times SO\left( 3,\mathbb{K}\right) .
\end{equation*}
\end{theorem}

\begin{proof}
Let $f\in \mathrm{Aut}(A\left( \frak{g}\right) )$ be given. Define $\hat{f}%
\in gl(\mathrm{Der}(A\left( \frak{g}\right) ))$ by $\hat{f}(T)=f\circ T\circ
f^{-1}$ for all $T\in \mathrm{Der}(A\left( \frak{g}\right) )$. Then $\hat{f}$
is an automorphism of $\mathrm{Der}(A\left( \frak{g}\right) )$. Now 
\begin{equation*}
\mathrm{Aut}(\mathrm{Der}(A\left( \frak{g}\right) ))=\mathrm{Aut}(\mathrm{%
diag}(\mathrm{ad}\left( \frak{g}\right) ))\times SO\left( 3,\mathbb{K}%
\right) ,
\end{equation*}
using Theorem \ref{thm-simple-der} and the fact that $\mathrm{Aut}(so\left(
3,\mathbb{K}\right) )=SO\left( 3,\mathbb{K}\right) $ \cite{jacobson}. We
have $\mathrm{diag}(\mathrm{ad}\left( \frak{g}\right) )\cong \mathrm{ad}%
\left( \frak{g}\right) $. Thus there exists $\hat{\phi}\in \mathrm{Aut}%
\left( \mathrm{ad}(\frak{g)}\right) $ such that 
\begin{equation*}
\hat{f}(\mathrm{diag}(\mathrm{ad}\ x))=\mathrm{diag}(\hat{\phi}(\mathrm{ad}\
x))
\end{equation*}
for all $x\in \frak{g}$. Since $\frak{g}$ is simple, $\frak{g}\cong \mathrm{%
ad}\left( \frak{g}\right) $. Thus define $\phi \in \mathrm{Aut}(\frak{g)}$
by $\hat{\phi}(\mathrm{ad}\ x)=\mathrm{ad}\ \phi (x)$ for all $x\in \frak{g}$%
. Then $\hat{\phi}(\mathrm{ad}\ x)=\phi \circ \mathrm{ad}\ x\circ \phi ^{-1}$
and hence 
\begin{equation*}
\mathrm{diag}(\hat{\phi}(\mathrm{ad}\ x))=\mathrm{diag}(\phi )\circ \mathrm{%
diag}(\mathrm{ad}\ x)\circ \mathrm{diag}(\phi ^{-1})
\end{equation*}
for all $x\in \frak{g}$. Next, there exists $R\in SO\left( 3,\mathbb{K}%
\right) $ such that 
\begin{equation*}
\hat{f}(M)=RMR^{-1}
\end{equation*}
for all $M\in so\left( 3,\mathbb{K}\right) $. Putting this together, we find
that 
\begin{equation*}
f\circ T\circ f^{-1}=\hat{f}(T)=(\mathrm{diag}(\phi )\circ R)\circ T\circ (%
\mathrm{diag}(\phi )\circ R)^{-1}
\end{equation*}
for all $T\in \mathrm{Der}(A\left( \frak{g}\right) )$, where $\mathrm{diag}%
(\phi )\circ R=R\circ \mathrm{diag}(\phi )$. By Lemma \ref{lemma-schur}, it
follows that $f^{-1}\circ \mathrm{diag}(\phi )\circ R=\lambda I$ for some $%
\lambda \in \mathbb{K}$. Since the left side is an automorphism, $\lambda =1$%
, and thus $f=\mathrm{diag}(\phi )\circ R$. This completes the proof.
\end{proof}

\end{document}